\documentclass{article}
\usepackage{amssymb, amsmath}
\usepackage[matrix, arrow, curve]{xy}
\newcommand \lap {\lambda^{\prime}}
\newcommand \pip {\pi^{\prime}}

\newcommand \la {\lambda}
\newcommand \gz {{\cal G}}
\newcommand \La {\Lambda}

\newcommand \R {{\cal R}}
\newcommand \A {{\cal A}}
\newcommand \W {{\cal W}}
\newcommand \T {{\cal T}}
\newcommand \q {{\bf q}}
\newcommand \Prob {\mathbb P}

\newcommand \U {{\cal U}}

\newcommand \Y {{\cal Y}}
\newcommand \sS {{\cal S}}

\newcommand \F {{\cal F}}
\newcommand \G {{\bf G}}
\newcommand \C {{\cal C}}

\newcommand \wab {{\cal W}_{{\cal A},B}}

\newcommand \cc {{\bf C}}
\newcommand \D  {{\cal D}} 
\newcommand \prt {{\bf p}_{r, \theta}}
\newcommand \Prt {{\bf P}_{r, \theta}}

\newcommand \jmtd {J^{(m)}_{\theta}\delta}
\newcommand \jmtpd {J^{(\pi^{-1}m)}_{-\theta}\delta}

\newcommand \aminus {{\bf a_{max}^{-}}}

\newcommand  {\modk} {{\cal  M}_{\kappa}}

\newcommand \HH {{\cal H}}

\newtheorem {lemma} {Lemma}
\newtheorem {theorem} {Theorem}
\newtheorem {proposition} {Proposition}
\newtheorem{corollary}{Corollary}
\begin{document}

\title{Decay of Correlations for the Rauzy--Veech--Zorich Induction Map 
on the Space of Interval Exchange Transformations and the 
Central Limit Theorem for the Teichm{\"u}ller Flow 
on the Moduli Space of Abelian Differentials.}

\author{Alexander I. Bufetov\footnote{Department of Mathematics, Princeton University}}
\maketitle

\section{Introduction}

The aim of this paper is to prove a stretched-exponential bound for the decay of correlations 
for the Rauzy-Veech-Zorich induction map on the space of interval exchange transformations 
(Theorem \ref{mainresult}). A Corollary is the Central Limit Theorem for the Teichm{\" u}ller flow
(Theorem \ref{clt-flow}).

The proof of Theorem \ref{mainresult} proceeds by  
approximating the induction map by a Markov chain satisfying the Doeblin condition, 
the method of Sinai \cite{sinai} and Bunimovich--Sinai \cite{bunimsinai}.
 The main ``loss of memory'' estimate is Lemma \ref{logsteps}.

\subsection{Interval exchange transformations.}

Let $m$ be a positive integer. Let $\pi$ be a permutation on $m$ symbols. 
The permutation $\pi$ will always be assumed {\it irreducible}, which means that 
$\pi\{1,\dots,k\}=\{1,\dots,k\}$ only if $k=m$. 
 
Let $\la$ be a vector in ${\mathbb R}^m_+$, $\la=(\la_1, \dots,\la_m)$, 
$\la_i>0$ for all $i$. Denote 
$$
|\la|=\sum_{i=1}^m\la_i.
$$

Consider the half-open interval $[0, |\la|)$.  Consider 
the points $\beta_i=\sum_{j<i}\la_j$, $\beta_i^{\pi}=\sum_{j<i}\la_{\pi^{-1}j}$.

Denote $I_i=[\beta_i, \beta_{i+1})$, $I_i^{\pi}=[\beta_i^{\pi}, \beta_{i+1}^{\pi})$.
The length of $I_i$ is $\la_i$, whereas the length of $I_i^{\pi}$  is $\la_{\pi^{-1}i}$. 

Set 

$$
T_{(\la,\pi)}(x)=x+\beta_{\pi i}^{\pi}-\beta_i {\rm \ for\ } x\in I_i.   
$$

The map $T_{(\la,\pi)}$ is called an interval exchange transformation corresponding to $(\la, \pi)$. 

The map $T_{(\la,\pi)}$ is an order-preserving isometry from $I_i$ onto $I_{\pi(i)}^{\pi}$.

We say that $\la$ is {\it irrational} if there are no
rational relations between $|\la|$, $\la_1$,$\la_2$, \dots $\la_{m-1}$. 

\begin{theorem}[Oseledets(\cite{oseledets}), Keane(\cite{keane})]

Let $\pi$ be irreducible and $\la$ irrational.
Then for any $x\in[0, \sum_{i=1}^{m}\la_i)$, the set $\{T_{(\la,\pi)}^n x, n\geq 0\}$ 
is dense in  $[0, \sum_{i=1}^{m}\la_i).$ 
\end{theorem}

\subsection{Rauzy operations $a$ and $b$.}
Let $(\la,\pi)$ be an interval exchange. Assume that $\pi$ is irreducible and  $\la$ is irrational.

Following Rauzy \cite{rauzy}, consider the induced map of $(\la,\pi)$ on the 
interval $[0, |\la|-min(\la_m, \la_{\pi^{-1}(m)}))$. The induced map is again 
an interval exchange of $m$ intervals. 
For $i,j=1,\dots,m$, 
denote by $E_{ij}$ an $m\times m$ matrix of which 
the $i,j$-th element is equal to $1$, all others to $0$.
Let $E$ be the $m\times m$-identity matrix.

\subsubsection{ Case $a$: $\la_{\pi^{-1}m}>\la_m$.} 

Define 
$$
A(a, \pi)=\sum_{i=1}^{\pi^{-1}(m)}E_{ii}+E_{m, \pi^{-1}m+1}+
\sum_{i=\pi^{-1}m+1}^m E_{i,i+1}
$$

$$
a\pi(j)=\begin{cases}  
\pi j,&\text{if $j\leq \pi^{-1}m$;}\\
                 \pi m,&\text{if $j=\pi^{-1}m+1$;}\\
                  \pi(j-1),&\text{ other $j$.}
\end{cases}
$$

If $\la_{\pi^{-1}m}>\la_m$, then the induced interval exchange 
of $T_{(\la, \pi)}$ on 
the interval
$[0, \sum_{i\neq m} \la_i)$ is $T_{(\lap, \pi^{\prime})}$, 
where $\lap=A(a,\pi)^{-1}\la$ and $\pip=a\pi$.

\subsubsection{ Case $b$: $\la_m>\la_{\pi^{-1}m}$.}

Define

$$
A(b, \pi)=E+E_{m,\pi^{-1}m}
$$

$$
b\pi(j)=\begin{cases}  
\pi j,&\text{if $\pi j\leq \pi m$;}\\
                 \pi j+1,&\text{if $\pi m<\pi j<m$;}\\
                  \pi m+1,&\text{ if $\pi j=m$.}
\end{cases}
$$

If $\la_{\pi^{-1}m}<\la_m$, then the induced interval exchange 
of $T_{(\la, \pi)}$ on 
the interval
$[0, \sum_{i\neq \Pi^{-1}m} \la_i)$ is $T_{(\lap, \pi^{\prime})}$, 
where $\lap=A(b,\pi)^{-1}\la$ and $\pip=b\pi$.

Note that  operations $a$ and $b$ are invertible on the space of permutations, namely, we have:

$$
a^{-1}\pi(j)=\begin{cases}
\pi(j),&\text{if $j\leq \pi^{-1}(m)$;}\\
\pi(j+1),&\text{if  $\pi^{-1}(m)+1<j<m$;}\\
\pi(\pi^{-1}(\pi(m)+1),&\text{if  $j=m$.}
\end{cases}
$$
$$
b^{-1}\pi(j)=\begin{cases}
\pi(j),&\text{ if $\pi(j)\leq \pi(m)$}\\
m,&\text{if $j=\pi^{-1}(\pi(m)+1)$;}\\
\pi(j)-1,&\text{if $\pi(j)>\pi(m)+1$.}
\end{cases}
$$

For $(\la,\pi)\in\Delta(\R)$, denote 

\begin{equation}
\label{TAB}
T_{a^{-1}}(\la,\pi)=(A(a^{-1}\pi, a)\la, a^{-1}\pi),\ 
T_{b^{-1}}(\la,\pi)=(A(b^{-1}\pi, b)\la, b^{-1}\pi).
\end{equation}

The interval exchange $T_{a^{-1}}(\la,\pi)$ is the preimage of $(\la,\pi)$ under the operation $a$, 
and the interval exchange $T_{b^{-1}}(\la,\pi)$ is the preimage of $(\la,\pi)$ under the operation $b$.

Normalize (dividing by $|\la|=\la_1+\dots +\la_m$) and set:

\begin{equation}
\label{tab}
t_{a^{-1}}(\la,\pi)=(\frac{A(a^{-1}\pi, a)\la}{|A( a^{-1}\pi, a)\la|}, a^{-1}\pi),\ 
t_{b^{-1}}(\la,\pi)=(\frac{A(b^{-1}\pi, b)\la}{|A( b^{-1}\pi, b)\la|}, b^{-1}\pi).
\end{equation}

\subsection{Rauzy class and Rauzy graph.}

If $\pi$ is an irreducible  permutation, then its 
{\it Rauzy class} is the set of all permutations
that can be obtained from $\pi$ by applying repeatedly 
the operations $a$ and $b$; the Rauzy class of the permutation 
$\pi$ is denoted ${\cal R}(\pi)$. 
Rauzy class has a natural 
structure of an oriented labelled graph: namely, the permutations 
of the Rauzy class are the vertices of the graph, and if 
$\pi= a\pi^{\prime}$ then we draw an edge from $\pi$ to $\pi^{\prime}$ 
and label it by $a$, and if $\pi=b\pi^{\prime}$ then we draw 
an edge from $\pi$ to $\pi^{\prime}$ and label it by $b$.
This labelled graph will be called the {\it Rauzy graph} of the permutation $\pi$.

For example, the Rauzy graph of the permutation $(4321)$ is 

$$
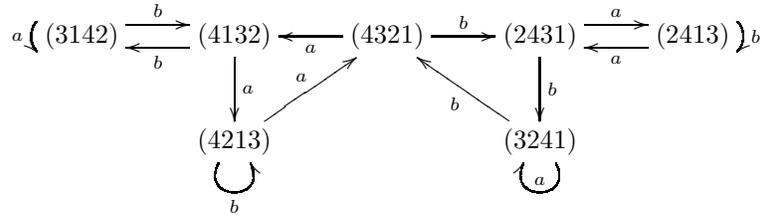
\xymatrix{
(3142)\ar@(ul,dl)[]_a\ar@<1ex>[r]^b&(4132)\ar@<1ex>[l]^b\ar[d]^a& 
(4321)\ar[l]^a\ar[r]^b&(2431)\ar[d]^b\ar@<1ex>[r]^a&(2413)\ar@<1ex>[l]^a 
\ar@(ur,dr)[]^b
\\
&(4213)\ar[ur]^a\ar@(dl,dr)[]_b&&(3241)\ar[ul]^b\ar@(dr,dl)[]_a& 
}
$$

For a permutation $\pi$, consider the set $\{a^n\pi, n\geq 0\}$. 
This set forms a cycle in the Rauzy graph which will be called 
the $a${\it -cycle} of $\pi$. Similarly, the set $\{b^n\pi, n\geq 0\}$
will be called the $b${\it-cycle} of $\pi$.

\subsection{The Rauzy-Veech-Zorich induction.}

Denote 
$$
\Delta_{m-1}=\{\la\in {\mathbb R}^m_+:|\la|=1\},
$$
$$
\Delta_{\pi}^+=\{\la\in\Delta_{m-1},\la_{\pi^{-1}m}>\la_m\},
\Delta_{\pi}^-=\{\la\in\Delta_{m-1},\la_m>\la_{\pi^{-1}m}\},
$$
$$
\Delta(\R)=\Delta_{m-1}\times {\cal R}(\pi).
$$

Define a map 
$$
{\cal T}: \Delta(\R)\to \Delta(\R)
$$

by 

$$
{\cal T}(\la,\pi)=\begin{cases}
(\frac{A(\pi, a)^{-1}\la}{|A(\pi, a)^{-1}\la|}, a\pi), &\text{if $\la\in\Delta_{\pi}^+$;}\\
(\frac{A(\pi, b)^{-1}\la}{|A(\pi, b)^{-1}\la|}, b\pi), &\text{if $\la\in\Delta_{\pi}^-$.}
\end{cases}
$$

Each $(\la, \pi)\in\Delta(\R)$ has exactly two preimages under the map ${\cal T}$, namely, 
$t_{a^{-1}}(\la,\pi)$ and $t_{b^{-1}}(\la,\pi)$ (\ref{tab}).

The set $\Delta(\R)$ is a finite union of simplices. 
Let ${\bf m}$ be the Lebesgue measure on $\Delta(\R)$ normalized in such a way
that ${\bf m}(\Delta(\R))=1$.

\begin{theorem}[Veech\cite{veech}]
The map ${\cal T}$ has an infinite  conservative ergodic invariant measure, absolutely continuous with respect to 
Lebesgue measure on $\Delta(\R)$. 
\end{theorem}

From this result Veech \cite{veech} derives that almost all (with respect to {\bf m}) interval exchange transformations 
are uniquely ergodic.

Denote 
$$
\Delta^+ = \cup_{\pip\in{\cal R}(\pi)}\Delta_{\pip}^+,
\Delta^- = \cup_{\pip\in{\cal R}(\pi)}\Delta_{\pip}^-.
$$

Following Zorich \cite{zorich}, we define the function $n(\la,\pi)$ in the following way.

$$
n(\la, \pi)=\begin{cases}  
\inf \{k>0:{\cal T}^k(\la,\pi)\in\Delta^-\},&\text{if $\la\in\Delta_{\pi}^+$;}\\
  \inf \{k>0: {\cal T}^k(\la,\pi)\in\Delta^+\},&\text{if $\la\in\Delta_{\pi}^-$.}
\end{cases}
$$ 

Define 
$$
{\cal G}(\la,\pi)={\cal T}^{n(\la,\pi)}(\la,\pi).
$$

The map $\gz$ will be referred to as {\it the Rauzy-Veech-Zorich induction map} \cite{rauzy, veech, zorich}.

For $(\la,\pi)\in\Delta(\R)$, denote 
$$
t_{a^{-n}}(\la,\pi)=t_{a^{-1}}^n(\la,\pi), t_{b^{-n}}(\la,\pi)=t_{b^{-1}}^n(\la,\pi), 
T_{a^{-n}}(\la,\pi)=T_{a^{-1}}^n(\la,\pi),T_{b^{-n}}(\la,\pi)=T_{b^{-1}}^n(\la,\pi).
$$

Under the map $\gz$, each interval exchange $(\la,\pi)$ has countably many preimages:  

$$
\gz^{-1}(\la, \pi)=\begin{cases}  
\{t_{a^{-n}}(\la,\pi), n\in{\mathbb N}\},&\text{if $(\la, \pi)\in\Delta^+$;}\\
                 \{t_{b^{-n}}(\la,\pi), n\in{\mathbb N}\},&\text{if $(\la, \pi)\in\Delta^-$.}
\end{cases}
$$

\begin{theorem}[Zorich\cite{zorich}]
The map ${\cal G}$ has an ergodic invariant probability measure, absolutely 
continuous with respect to Lebesgue on $\Delta({\cal R})$.   
\end{theorem}

Denote this invariant measure by $\nu$; the probability with respect to $\nu$ will be denoted by $\Prob$.

Let $\rho(\la, \pi)$ be the density of $\nu$ with 
respect to the Lebesgue measure ${\bf m}$. Zorich \cite{zorich} showed that for any $\pi\in{\cal R}$ 
there exist two positive rational homogeneous of degree $-m$  functions $\rho^+_{\pi}$, 
$\rho^{-}_{\pi}$ such that 

\begin{equation}
\label{rhopm}
\rho(\la, \pi)=\begin{cases}  
\rho^+_{\pi}(\la),&\text{if $\la\in\Delta_{\pi}^+$;}\\
                 \rho^{-}_{\pi}(\la),&\text{if $\la\in\Delta_{\pi}^-$.}
\end{cases}
\end{equation}

{\bf Remark.} In particular, the invariant density is bounded from below: there exists 
a positive constant $C(\R)$, depending on the Rauzy class only and such that 
$\rho(\la,\pi)\geq C(\R)$ for any $(\la,\pi)\in\Delta(\R)$.

The map $\gz$ is not mixing: indeed, from the definition of $\gz$, we have 
$$
\gz(\Delta^+)=\Delta^-,\  
\gz(\Delta^-)=\Delta^+.
$$

Let ${\cal B}$ be the Borel $\sigma$-algebra on $\Delta(\R)$, and let ${\cal B}_n=\gz^{-n}{\cal B}$.
We have  ${\cal B}_{n+2}\subset {\cal B}_n$.
Recall \cite{rokhlin} that {\it exactness} of the map $\gz^2$ means, by definition, that 
the $\sigma$-algebra $\cap_{n=1}^{\infty} {\cal B}_{2n}$ 
is trivial \cite{rokhlin} (in other words, that Kolmogorov's $0-1$ law holds for the map $\gz^2$.)

\begin{proposition} 
The map ${\cal G}^2:\Delta^+\to \Delta^+$ is exact with respect to $\nu|_{\Delta^+}$. 
\end{proposition}

This Proposition is proven in Section \ref{proofexactness}; it implies strong mixing for the map $\gz^2$.

\subsection{The main result} 

Introduce a metric on $\Delta_{m-1}$ by setting
\begin{equation}
\label{birkhoffmetric}
d(\la, \lap)=\log \frac{\max_i\frac{\la_i}{\lap_i}}{\min_i\frac{\la_i}{\lap_i}}.
\end{equation}

Now introduce a metric on $\Delta(\R)$ by setting 
$$
d((\la,\pi), (\lap, \pip))=\begin{cases}  
2+d(\la,\lap),&\text{if $\pi\neq\pip$;}\\
                 d(\la,\lap),&\text{if $\pi=\pip$.}
\end{cases}
$$

For $\alpha>0$, let $H_{\alpha}$ be the space of functions 
$\phi:\Delta(\R)\to{\mathbb R}$ such that if $d((\la,\pi), (\lap, \pip)\leq 1$, then
$|\phi(\la,\pi)-\phi(\lap, \pip)|\leq C d((\la,\pi),(\lap, \pip))^{\alpha}$ for some constant $C$.

Define 
$$
C_{H_{\alpha}}(\phi)=\max_{d((\la,\pi), (\lap, \pip))\leq 1}
\frac{|\phi(\la,\pi)-\phi(\lap, \pip)|} {d((\la,\pi), (\lap, \pip))^{\alpha}}, 
$$

The main result of this paper is

\begin{theorem}
\label{mainresult}
Let $\gz:\Delta(\R)\to\Delta(\R)$ be the Rauzy-Veech-Zorich 
induction map and let $\nu$ be the absolutely continuous  
invariant measure.  

Let $p>2$.
Then, for any $\alpha>0$, there exist positive constants $C, \delta$ such that 
for any  $\phi\in H_{\alpha}\cap L_p(\Delta^+(\R), \nu)$ and $\psi\in L_2(\Delta^+(\R), \nu)$ we have 

$$
|\int \phi\ \times \psi\circ \gz^{2n} d\nu-\int\phi d\nu\int\psi d\nu|\leq C\exp(-\delta n^{1/6})
(C_{H_{\alpha}}(\phi)+|\phi|_{L_p})(|\psi|_{L_2}).  
$$
\end{theorem}

Denote by ${\cal N}(0, \sigma)$ the Gaussian distribution with mean $0$ and variance $\sigma$.
By \cite{gordin, liverani, viana}, we have 

\begin{corollary}
\label{clt-discrete}

Let $\phi\in H_{\alpha}\cap L_p(\Delta(\R)^+, \nu)$, $\int \phi d\nu=0$. 
Assume that there does not exist $\psi\in L_2(\Delta(\R)^+, \nu)$
such that $\phi=\psi\circ \gz^2-\psi$. Then there exists $\sigma>0$ such that 

$$
\frac 1 {\sqrt N} \sum_{n=0}^{N-1} \phi\circ \gz^{2n} \xrightarrow{d} {\cal N}(0, \sigma) \ {\rm as }\  N\to\infty.
$$
\end{corollary}

\subsection{Veech's space of zippered rectangles}

A zippered rectangle associated to the Rauzy class $\R$ is a quadruple $(\la, h, a, \pi)$,  where $\la \in 
{\mathbb R}_+^m, 
h\in{\mathbb R}^m_+, a\in{\mathbb R}^m, \pi\in{\cal R}$, and the vectors 
$h$ and $a$ satisfy the following equations and inequalities (one introduces 
auxiliary components $a_0=h_0=a_{m+1}=h_{m+1}=0$, and sets $\pi(0)=0$, 
$\pi^{-1}(m+1)=m+1$.): 

$$
h_i-a_i=h_{\pi^{-1}(\pi(i)+1)}-a_{\pi^{-1}(\pi(i)+1)-1}, i=0, \dots, m
$$
$$
h_i\geq 0,   i=1, \dots,m, \  
a_i\geq 0,  i=1, \dots, m-1,  
$$
$$
a_i\leq \min (h_i, h_{i+1}) {\rm \ for\ } i\neq m,i\neq \pi^{-1}m,  
$$
$$
a_m\leq h_m, \ a_m\geq -h_{\pi^{-1}m},\ 
a_{\pi^{-1}m}\leq h_{\pi^{-1}m+1} 
$$

The area of a zippered rectangle is given by the expression $\la_1h_1+\dots+\la_mh_m$.
Following Veech, we denote by $\Omega(\R)$ the space of all zippered 
rectangles, corresponding to a given Rauzy class $\R$ and satisfying the condition 
$$
\la_1h_1+\dots+\la_mh_m=1.
$$
We shall denote 
by $x$ an individual zippered rectangle.

Veech further defines a map $\U$ and a flow $P^t$ on the space of 
zippered rectangles in the following way:

$$
P^t(\la,h,a,\pi)=(e^t\la, e^{-t}h, e^{-t}a, \pi).
$$

$$
\U(\la,h,a,\pi)=\begin{cases}  
(A^{-1}(a,\pi)\la, A^t(a,\pi)h, a^{\prime}, a\pi), &\text{if $(\la,\pi)\in\Delta^-$}\\
 (A^{-1}(b,\pi)\la, A^t(b,\pi)h, a^{\prime\prime}, b\pi),&\text{if $(\la,\pi)\in\Delta^+$},
\end{cases}
$$

where 

$$
a^{\prime}_i=\begin{cases}  
a_i, &\text{if $j<\pi^{-1}m$,}\\
h_{\pi^{-1}m}+a_{m-1}, &\text {if $i=\pi^{-1}m$,}\\
a_{i-1}, &\text{other $i$}.
\end{cases}
$$

$$
a^{\prime\prime}_i=
\begin{cases}  
a_i, &\text{if $j<m$,}\\
-h_{\pi^{-1}m}+a_{\pi^{-1}m-1}, 
&\text {if $i=m$.}
\end{cases}
$$

The map $\U$ is invertible; $\U$ and $P^t$ commute (\cite{veech}). 

Denote

$$
\tau(\lambda, \pi)=(\log(|\la|-\min(\la_m,\la_{\pi^{-1}m})),
$$

and for $x\in \Omega(\R)$, $x=(\la, h, a ,\pi)$, write
 
$$
\tau(x)=\tau(\la,\pi).
$$

Now define 

$$
\Y(\R)=\{x\in\Omega(\R):|\la|=1\}.
$$

and 
$$
\Omega_0(\R)=\bigcup_{x\in\Y(\R), 0\leq t\leq \tau(x)}P^tx.
$$
$\Omega_0(\R)$ is a fundamental domain for $\U$ and,
identifying the points $x$ and $\U x$ in $\Omega_0(\R)$, we obtain 
a natural flow, also denoted by $P^t$, on $\Omega_0(\R)$.

The space $\Omega(\R)$ has a natural Lebesgue measure class and so does  
the transversal $\Y(\R)$. 
Veech \cite{veech} has proved the following Theorem.

\begin{theorem}
\label{zipmeas}
There exists a measure $\mu_{\R}$ on $\Omega(\R)$, absolutely continuous with respect to Lebesgue, preserved by both the map 
$\U$  and the flow $P^t$ and such that $\mu_{\R}(\Omega_0(\R))<\infty$. 
\end{theorem}

For $x\in\Y(\R)$, define 

$$
\sS(x)=\U P^{\tau(x)}(x).
$$

The map $\sS$ is a lift of ${\cal T}$ to the space of zippered rectangles: 
indeed, if 

$$
\sS(\la,h,a,\pi)=(\lap, h^{\prime}, a^{\prime}, \pip),
$$

then $(\lap, \pip)={\cal T}(\lap, \pip)$.

Since $\Y(\R)$ is a transversal to the flow, the measure $\mu_{\R}$ induces an absolutely 
continuous measure $\mu_{\R}^{(1)}$ on $\Y(\R)$; since $\mu_{\R}$ is both $\U$ and $P^t$-invariant, the 
measure $\mu_{\R}^{(1)}$ is $\sS$-invariant. Since $\mu_{\R}(\Omega_0(\R))<\infty$, the measure $\mu_{\R}^{(1)}$
is conservative; it is, however, infinite (Veech \cite{veech}).

Zorich \cite{zorich} constructed a different section for the flow $P^t$, for which the restricted measure has finite total mass. 

Following Zorich \cite{zorich},  define 
$$
\Omega^+(\R)=\{ x=(\la, h,a,\pi):(\la,\pi)\in\Delta^+, a_m\geq 0 \}.
$$

$$
\Omega^-(\R)=\{ x=(\la, h,a,\pi):(\la,\pi)\in\Delta^-, a_m\leq 0 \},
$$
$$
\Y^+(\R)=\Y(\R)\cap \Omega^+(\R),\  \Y^-(\R)=\Y(\R)\cap \Omega^-(\R),\  \Y^{\pm}(\R)=\Y^+(\R)\cup \Y^-(\R).
$$

Take $x\in\Y^{\pm}(\R)$,  $x=(\la,h,a,\pi)$, and define  

$$
\F(x)=\sS^{n(\la,\pi)}x.
$$

The map $\F$ is a lift of the map $\gz$ to the space of zippered rectangles:
if 

$$
\F(\la,h,a,\pi)=(\lap, h^{\prime}, a^{\prime}, \pip),
$$

then $(\lap, \pip)=\gz(\lap, \pip)$.

We shall see, moreover, that the map $\F$ can be almost surely (with respect to Lebesgue)
identified with the natural extension of the map $\gz$ (Section 3).

If $x\in\Y^+$, then $\F(x)\in \Y^-$, and if $x\in\Y^-$, then $\F(x)\in 
\Y^+$.
The map $\F$ is the induced map of $\sS$ to the subset 
$\Y^{\pm}(\R)$.

Since $\Y^{\pm}(\R)$ is a transversal to the flow $P^t$, the measure 
$\mu_{\R}$ naturally induces an absolutely 
continuous measure ${\overline \nu}$ on $\Y^{\pm}(\R)$; since $\mu_{\R}$ is both $\U$ and $P^t$-invariant, the 
measure ${\overline \nu}$ is $\F$-invariant.

Zorich \cite{zorich} proved

\begin{theorem}
The measure ${\overline \nu}$ is finite and ergodic for $\F$.
\end{theorem}

Since the map $\gz$ is exact (as is shown in Section 4), 
the map $\F$ satisfies the $K$-property of Kolmogorov, and, in particular, is strongly mixing. 
Decay of correlations is proven for the map $\F$ as well. 

Introduce a metric on the space of zippered rectangles in the following way.
Take two zippered rectangles 
 $x=(\la,h,a, \pi)$ and $x^{\prime}=(\lap, h^{\prime}, a^{\prime}, \pip)$. 
Write 

$$
d((\la,h,a), (\lap, h^{\prime}, a^{\prime}))=
\log \frac{\max_i\frac{\la_i}{\lap_i},\frac {h_i}{h^{\prime}_i}, 
\frac{|a_i|}{|a^{\prime}_i|}, \frac{|h_i-a_i|}{|h^{\prime}_i-a^{\prime}_i|}}
{\min_i\frac{\la_i}{\lap_i},\frac {h_i}{h^{\prime}_i}, 
\frac{|a_i|}{|a^{\prime}_i|}, \frac{|h_i-a_i|}{|h^{\prime}_i-a^{\prime}_i|}}.
$$

Define the metric on $\Omega(\R)$ by 

$$
d(x,x^{\prime})=\begin{cases}  
d((\la,h, a), (\lap, h^{\prime}, a^{\prime})
&\text{if $\pi=\pip$ and $\frac {a_m}{a^{\prime}_m}>0$;}\\
2+d((\la,h,a), \lap, h^{\prime}, a^{\prime}), &\text{otherwise}.
            
\end{cases}
$$

As above, for $\alpha>0$, let $H_{\alpha}$ be the space of functions 
$\phi:\Y^{\pm}(\R)\to{\mathbb R}$ such that if $d(x, x^{\prime})\leq 1$, then
$|\phi(x)-\phi(x^{\prime})|\leq C d(x,x^{\prime})^{\alpha}$ for some constant $C$.

Note that the distance  $d(x,x^{\prime})$ is not defined if $a_i=0$ or $a^{\prime}_i=0$ 
for some $i=1, \dots, m$; nothing, therefore, is said about the values of a function from 
$H_{\alpha}$ at such points. This does not represent a problem, however, since 
we only need the space $H_{\alpha}$ for the Central Limit Theorem, and 
for for such a result we may deal with functions defined almost everywhere.

Define 
$$
C_{H_{\alpha}}(\phi)=\max_{d(x, x^{\prime})\leq 1}
\frac{|\phi(x)-\phi(x^{\prime})|} {d(x,x^{\prime})^{\alpha}}.
$$

\begin{theorem}
\label{zipdecay}
Let $\F:\Y^{\pm}(\R)\to\Y^{\pm}(\R)$ be 
the Rauzy-Veech-Zorich 
induction map on the space of zippered rectangles
and let ${\overline \nu}_{\R}$ be the absolutely continuous  
invariant probability measure.  
Let $p>2$. 
Then, for any $\alpha>0$, there exist positive constants $C, \delta$ such that 
for any  $\phi, \psi\in H_{\alpha}\cap L_p(\Y(\R), {\overline \nu}_{\R})$ we have 

$$
|\int \phi\ \times \psi\circ \F^{2n} d{\overline \nu}_{\R}-
\int\phi d{\overline \nu}_{\R}\int\psi d{\overline \nu}_{\R}|\leq C\exp(-\delta n^{1/6})
(C_{H_{\alpha}}(\phi)+|\phi|_{L_p})(C_{H_{\alpha}}(\psi)+|\psi|_{L_p})  
$$
\end{theorem}

Theorem \ref{zipdecay} will be established simultaneosuly with the Theorem \ref{mainresult}. 
Indeed, the map $\F$ can be almost surely identified with the natural extension of the map $\gz$, and
the method of Markov approximations of 
of Sinai \cite{sinai} and Bunimovich--Sinai \cite{bunimsinai} allows to 
obtain the decay of correlations for the invertible case simultaneously with that for 
the noninvertible one.

Since the flow $P^t$ is a special flow over the map $\F$, by 
the Theorem of Melbourne and T{\"o}r{\"o}k \cite{torok}, the decay of correlations for the map $\F$ 
allows to obtain the Central Limit Theorem for the flow $P^t$.  

Denote by $X_t$ the derivative with respect to the flow $P^t$.

\begin{theorem}
\label{cltpt}

Let $p>2$ and let $\phi\in H_{\alpha}(\Omega_0(\R))\cap L_p(\Omega_0(\R), \mu_{\R})$ 
satisfy $\int \phi d\nu=0$. 
Assume that there does not exist  $\psi\in L_2(\Omega_0(\R), \mu_{\R})$
such that $\phi=X_t\psi$. Then there exists $\sigma>0$ such that 

$$
\frac 1 {\sqrt T} \int_0^T \phi\circ P^t \xrightarrow{d} {\cal N}(0, \sigma) \ {\rm as } \ T\to\infty.  
$$
\end{theorem}

This Theorem will be proved in Section 16.

\subsection{Zippered rectangles and the moduli space of holomorphic 
differentials.}

Let $g\geq 2$ be an integer. Take an arbitrary integer vector $\kappa=(k_1, \dots, k_{\sigma})$ 
such that $k_i>0$, $k_1+\dots +k_{\sigma}=2g-2$. 

Denote by $\modk $ the moduli space 
of Riemann surfaces of genus $g$ endowed with a holomorphic differential  
of area $1$ with singularities of orders $k_1, \dots, k_{\sigma}$. 
(the {\it stratum} in the moduli space of holomorphic differentials).
Denote by $g_t$ the Teichm{\"u}ller flow on $\modk$ 
(see \cite{forni}, \cite{masur}, \cite{hubmas}, \cite{kz1}).
The flow $g_t$ preserves a natural absolutely continuous probability
measure on $\modk$ (\cite{masur},\cite{veech}, \cite{kz1}). We denote that measure by $\mu_{\kappa}$.

A zippered rectangle naturally defines a Riemann surface endowed with 
a holomorphic differential of area $1$. The orders of the 
singularities of $\omega$ are uniquely 
defined  by the Rauzy class of the permutation $\pi$ (\cite{veech}).

For any $\R$ we thus have a map 

$$
\pi_{\R}: \Omega_{\R}\rightarrow\modk,
$$

where ${\kappa}$ is uniquely defined by $\R$.

Veech \cite{veech} proved

\begin{theorem}[Veech]
\label{zipmodule}
\begin{enumerate}
\item The set $\pi_0(\Omega_0(\R))$ is a connected component of $\modk$. 
Any connected component of any $\modk$ has the form $\pi_0(\Omega_0(\R))$ 
for some $\R$.
\item The map $\pi_0$ is finite-to-one and almost everywhere locally 
bijective.
\item $\pi_0(\U x)=\pi_0(x)$.
\item The flow $P^t$ on $\Omega_0(\R)$ projects under $\pi_0$ 
to the Teichm{\"u}ller flow $g_t$ on the corresponding connected component of 
$\modk$.
\item $(\pi_{\R})_*\mu_{\kappa}=\mu_{\R}$.
\end{enumerate}
\end{theorem}

A detailed treatment of the relationship between Rauzy classes, zippered rectangles and connected components 
is given by M.Kontsevich and A.Zorich in \cite{kontsevichzorich}.

Say that a function $\psi:\modk\to{\mathbb R}$ is {\it H{\"o}lder in the sense of Veech}
if there exists a H{\"o}lder function $\phi:\Omega_0(\R)\to {\mathbb R}$ such that 
$\psi\circ\pi_0=\phi$. 

{\bf Remark.} This definition has a natural interpretation in 
terms of cohomological coordinates of Hubbard and Masur \cite{hubmas}. Indeed, under the map $\pi_0$ 
the Veech coordinates on the space of zippered rectangles
correspond, upto a linear change of variables, to 
the cohomological coordinates of Hubbard and Masur.
Locally, one can associate a Hilbert metric to those coordinates. A function H{\"o}lder in the sense 
of Veech if and only if it is H{\"o}lder with respect to that metric.
Note that the thus defined local Hilbert distance between two elements in $\modk$
majorates the Teichm{\"u}ller distance between their underlying surfaces.
Therefore, if a function $\phi: \modk\to{\mathbb R}$ is a lift of a smooth 
function from the underlying moduli space 
${\cal M}_g$ of compact surfaces of genus $g$, then $\phi$ is H{\"o}lder in the sense of Veech.

Denote by ${\cal X}_t$ the derivative in the direction of the flow $g_t$. 

Theorem \ref{cltpt} and Theorem \ref{zipmodule} imply the following

\begin{theorem}
\label{clt-flow}
Let $\HH$ be a connected component of $\modk$.
Let  $p>2$, and let $\psi\in L_p(\HH, \mu_{\kappa})$ be  H{\"o}lder in the sense of Veech and satisfy
$\int \phi d\mu_{\kappa}=0$. 
Assume that there does not exist  $\psi\in L_2(\HH, \mu_{\kappa})$
such that $\phi={\cal X}_t\psi$. Then there exists $\sigma>0$ such that 

$$
\frac 1 {\sqrt T} \int_0^T \phi\circ g_t dt \xrightarrow{d} {\cal N}(0, \sigma) \ {\rm as } \ T\to\infty.  
$$

\end{theorem}

\subsection{Outline of the Proof of Theorem \ref{mainresult}.}

First, one takes a subset of the 
space $\Delta(\R)$ such that the induced map of ${\cal G}$ is uniformly expanding
(namely, the set of all interval exchanges such that the renormalization matrix 
for them is a fixed matrix all whose elements are positive, see Proposition \ref{unifexp}; note that the return map on such a subset  
is an essential element in Veech's proof of unique ergodicity \cite{veech}). 
Then one  estimates the statistics of return times in this subset, in the spirit of 
Lai-Sang Young \cite{lsy}. 
After that, the method of Markov approximations, due to Sinai \cite{sinai}, Bunimovich 
and Sinai \cite{bunimsinai}, is used to complete the proof.

The paper is organized as follows. In Section 2, we state auxiliary propositions about unimodular 
matrices. In Section 3, following Veech \cite{veech} and Zorich \cite{zorich}, we construct 
symbolic dynamics for the Rauzy-Veech-Zorich induction map $\gz$,  
compute its transition probabilities in the sense of Sinai \cite{sinai}, and identify 
the natural extension of $\gz$ with $\F$.
In Section 4, we establish the exactness of $\gz^2$. 
In Section 6, we state the main
Lemma \ref{logsteps}, whose proof takes  
Sections 6 -- 10.
In the remainder of the paper we apply the Markov approximation method of Sinai \cite{sinai},
Bunimovich and Sinai \cite{bunimsinai}, in order to obtain the decay of correlations for $\gz$ and $\F$. 
In the final Section, we apply the Theorem of Melbourne and T{\"o}r{\"o}k to obtain the Central Limit Theorem 
for the Teichm{\"u}ller flow.

\section{Matrices}

Let $A$ be an $m\times m$-matrix with positive entries.

Denote 

$$
|A|=\sum_{i,j=1}^m A_{ij}
$$

$$
col(A)=\max_{i,j,k} \frac{A_{ij}}{A_{kj}},
$$

$$
row(A)=\max_{i,j,k}\frac{A_{ij}}{A_{ik}}
$$

\begin{proposition}
Let $Q$ be a matrix with positive entries, $A$ a matrix with nonnegative entries without zero 
columns or rows.

Then all entries of the matrices $AQ$ and $QA$ are positive, and, moreover, we have

$$
row(AQ)\leq row(Q), col(QA)\leq col(Q)
$$
\end{proposition}

\begin{corollary}
Let $Q$ be a matrix with positive entries, $A$ a matrix with nonnegative entries without zero 
columns or rows.
$$
row(QAQ)\leq row(Q), col(QAQ)\leq col(Q)
$$
\end{corollary}

Let $A$ be an $m\times m$ matrix with nonnegative entries and determinant $1$. Consider the map 
$J_A:\Delta_{m-1}\to\Delta_{m-1}$ given by 

$$
J_A(\la)=\frac{A\la}{|A\la|}.
$$

Then

\begin{equation}
\label{detja}
det DJ_A(\la)=\frac 1{|A\la|^m}.
\end{equation}

Suppose all entries of $A$ are positive; then, for any $\la,\lap\in\Delta_{m-1}$, we have
\begin{equation}
\label{distortion}
row(A)^{-m}\leq\frac{det DJ_A(\la)}{det DJ_A(\lap)}\leq row(A)^m,
\end{equation}

whence we have the following 

\begin{proposition}
\label{bdddist}
Let $C\subset\Delta_{m-1}$ and let $A$ be a matrix with positive entries and determinant $1$.
Then 
$$
row(A)^{-m}\frac{{\bf m}(C_1)}{{\bf m}(C_2)}
\leq \frac{{\bf m}(J_A(C_1))}{{\bf m}(J_A(C_2))}\leq 
row(A)^m\frac{{\bf m}(C_1)}{{\bf m}(C_2)}.
$$
\end{proposition}

We also note the following well-known Lemma (see, for example, \cite{viana}):

\begin{lemma}
\label{unifcontract}
Suppose all entries of the matrix $A$ are positive. 
Then the map $J_A$ is uniformly contracting with respect to the Hilbert metric.
\end{lemma}

\section{Symbolic dynamics for $\gz$.}

First, following Veech \cite{veech} and Zorich \cite{zorich}, we describe a Markov partition and a
symbolic dynamics for the map $\gz^2$, then we identify almost surely 
the induction map $\F$ on the space of zippered rectangles with the natural extension of 
$\gz$, and, finally,  we  compute for $\gz$ its transition probabilities in
the sense of Sinai \cite{sinailectures}.

\subsection{The alphabet}

Let $\pi\in {\cal R}$, and let $n$ be a positive integer. 

Set $$\La(a, n, \pi)=\{\la: {\rm there \ exists} \ 
(\lap, \pi^{ \prime})\ \ {\rm such \ that} \ 
\lap\in\Delta_{\pi^{\prime}}^+ \ {\rm and}\  
(\la, \pi)=t_{a^{-n}}(\lap, \pip)\}$$

$$
\Delta(a,n,\pi)=\{(\la,\pi), \la\in\La(a,n, \pi)\}
$$
 
In other words, $\Delta(a,n,\pi)$ is the set of interval exchange transformations 
such that the application of the Zorich induction results in the 
application of the $a$-operation $n$ times.

The sets $\Delta(a,n,\pi)$ and $\Delta(a, n^{\prime},  \pip)$ are disjoint
unless $n=n^{\prime}$, $\pi=\pip$, and
 
$$
\Delta_{\pi}^-=\cup_{n=1}^{\infty} \Delta(a,n, \pi) 
$$
up to a set of measure zero (namely, a union of countably many hyperplanes 
on which  Zorich induction is not defined).

If $\pip=a^n\pi$, then we have 
$$
\gz \Delta(a,n, \pi)=\Delta_{\pip}^+.
$$ 

Similarly,  for $\pi\in {\cal R}$, and  $n$ a positive integer,
set 

$$\La(b,n, \pi)=\{\la: {\rm there \ exists}\  
(\lap, \pi^{ \prime})\  {\rm such \ that}\   
\lap\in\Delta_{\pi^{\prime}}^- \ {\rm and}\  
(\la, \pi)=t_{b^{-n}}(\lap, \pip)\}.
$$

$$
\Delta(b,n, \pi)=\{(\la,\pi), \la\in\La(b,n, \pi)\}.
$$
 
In other words, $\Delta(b,n,\pi)$ is the set of interval exchange transformations 
such that the application of the Zorich induction results in the 
application of the $b$-operation $n$ times.

The sets $\Delta(b,n, \pi)$ and $\Delta(b, n^{\prime},  \pip)$ are disjoint
unless $n=n^{\prime}$, $\pi=\pip$, and
 
$$
\Delta_{\pi}^+=\cup_{n=1}^{\infty} \Delta(b,n, \pi)
$$

up to a set of measure zero (namely, a union of countably many hyperplanes
on which the Zorich induction is not defined).

If $\pip=b^n\pi$, then, clearly, 
$$
\gz (\Delta(b,n, \pi))=\Delta_{\pip}^-.
$$ 

Note that the sets $\Delta(a,n,\pi)$ and $\Delta(b, n^{\prime},\pip)$ are always disjoint,
since we have $\Delta(a,n,\pi)\subset \Delta_{\pi}^-$, 
$\Delta(b, n^{\prime},\pip)\subset\Delta_{\pip}^+$.

The sets $\Delta(a,n, \pi)$, $\Delta(b,n, \pi)$, for all  $n>0$ and 
all $\pi\in{\cal R}$, form a Markov partition for $\gz$.

\subsection{Words}

Consider the alphabet

$$
{\cal A}=\{ (c,n,\pi), c=a \ {\rm or}\  b \}
$$

For $w_1\in\A$, $w_1=(c_1, n_1, \pi_1)$, we write $c_1=c(w_1), \pi_1=\pi(w_1), n_1=n(w_1)$.

For $w_1,w_2\in\A$, $w_1=(c_1, n_1, \pi_1)$, $w_2=(c_2, n_2, \pi_2)$, define 
the function $B(w_1, w_2)$ in the following way:
$B(w_1, w_2)=1$ if $c_1^{n_1}\pi_1=\pi_2$ and $c_1\neq c_2$ and $B(w_1, w_2)=0$ otherwise.

Let 

$$
W_{{\cal A}, B}=\{w=w_1\dots w_n , w_i\in{\cal A}, B(w_i, w_{i+1})=1 \  {\rm for \ all}\  i=1, \dots,n\}.
$$

For $w_1\in\A$, $w_1=(c_1, n_1, \pi_1)$, set 
$$
A(w)=A(c_1, c_1^{-n_1}\pi_1)\dots A(c_1, c_1^{-1}\pi_1)A(c_1, \pi_1),
$$
and for $w\in W_{\A,B}$, $w=w_1\dots w_n$, set 
$$
A(w)=A(w_1)\dots A(w_n).
$$

Also, for $w_1\in \A$, $\pi\in\R$, set $w_1^{-1}\pi=c_1^{-n_1}\pi$, and for 
$w\in W_{\A,B}$, $w=w_1\dots w_n$, set 
$$
w^{-1}\pi=w_1^{-1}\dots w_n^{-1}\pi.
$$

For $w\in W_{\A, B}$, define a map $t_w:\Delta(\R)\to\Delta(\R)$ by 

$$
t_w(\la,\pi)=(\frac{A(w)\la}{|A(w)\la|}, w^{-1}\pi)
$$

Consider also the map

$$
T_w(\la,\pi)=({A(w)\la}, w^{-1}\pi)
$$

For $w_1\in\A$, $w_1=(c_1, n_1, \pi_1)$, we write $\Delta(w_1)=\Delta(c_1, n_1\pi_1)$.

For $w\in W_{\A,B}$, $w=w_1\dots w_n$, denote 

$$
\Delta(w)=t_w(\Delta(\R)).
$$

Then, by definition, 

$$
\Delta(w)=\{(\la,\pi): (\la,\pi)\in\Delta(w_1), \gz(\la,\pi)\in\Delta(w_2), \dots, \gz^{n-1}(\la,\pi)\in\Delta(w_n)\}.
$$

Say that $w_1\in\A$ is compatible with $(\la,\pi)\in\Delta(\R)$ if 

\begin{enumerate}
\item either $\la\in\Delta_{\pi}^{+}$, $c_1=a$, and $a^{n_1}\pi_1=\pi$
\item or $\la\in\Delta_{\pi}^{-}$, $c_1=b$, and $b^{n_1}\pi_1=\pi$.
\end{enumerate}

Say that a word $w\in W_{\A,B}$, $w=w_1\dots w_n$ is compatible with $(\la,\pi)$ if 
$w_n$ is compatible with $(\la,\pi)$.  

We can write 

$$
\gz^{-n}(\la,\pi)=
\{t_w(\la,\pi): |w|=n \ {\rm and}\  w\  {\rm is\  compatible \ with\ } \ (\la,\pi)\}.
$$

Suppose that a word $w\in  W_{\A,B}$ is compatible with both $(\la,\pi)$ and 
$(\lap, \pi)$. Then

$$
d(t_w(\la,\pi), t_w(\lap, \pi))\leq d((\la,\pi), (\lap, \pip)).
$$ 

If, moreover, all entries of the the matrix $A(w)$ are positive, then, by Lemma \ref{unifcontract}, there 
exists $\alpha(w)$, $0<\alpha(w)<1$, such that 

$$
d(t_w(\la,\pi), t_w(\lap, \pi))\leq \alpha(w)d((\la,\pi), (\lap, \pip)).
$$

We therefore have 

\begin{proposition}
\label{unifexp}
Let $w\in  W_{\A,B}$ be such that all entries of the matrix $A(w)$ are positive.
Then the return map of $\gz$ on $\Delta(w)$ is uniformly expanding with respect to the 
Hilbert metric.
\end{proposition}

\subsection{Sequences} 

Now let $$\Omega_{\A, B}=\{\omega=\omega_1\dots \omega_n\dots ,\  \omega_n\in\A, 
B(\omega_n, \omega_{n+1})=1  \ {\rm for \ all}\  n\in {\mathbb N}\}$$

and 

$$\Omega_{\A, B}^{\mathbb Z}=\{\omega=\dots \omega_{-n}\dots \omega_1\dots \omega_n\dots ,\  \omega_n\in\A, 
B(\omega_n, \omega_{n+1})=1 
\ {\rm for \ all}\  n\in{\mathbb Z}\}$$

Denote by  $\sigma$ the shift on both these spaces.

There is a natural map $\Phi:\Delta\to\Omega_{\A, B}$ given by the formula 

$$
\Phi(\la,\pi)=\omega_1 \dots \omega_n \dots 
$$

if 

$$
\gz^n(\la,\pi)\in\Delta(\omega_n)
$$

The measure $\nu$ projects under $\Phi$ to a $\sigma$-invariant measure on $\Omega_{\A, B}$; probability with respect to 
that measure will be denoted by $\Prob$.

For $w\in W_{\A,B}$, $w=w_1\dots w_n$, let 
$$
C(w)=\{\omega\in\Omega_{\A,B}:\omega_1=w_1, \dots, \omega_n=w_n\}.
$$
We have then 
$$
\Delta(w)=\Phi^{-1}(C(w)).
$$

W. Veech \cite{veech} has proved the following

\begin{proposition}
The map $\Phi$ is $\nu$-almost surely bijective.
\end{proposition}

We thus obtain a symbolic dynamics for the map $\gz$.

\subsection{ The natural extension.}

Consider the natural extension for the map $\gz$. 

The phase space is the space of sequences of interval exchanges; it will be convenient to number them 
by negative integers. We set: 

$$
{\overline \Delta}(\R)=
$$
$$\{{\bf x}=\dots (\la(-n), \pi(-n)), \dots, (\la(0), \pi(0))| \ 
\gz(\la(-n), \pi(-n))=(\la(1-n), \pi(1-n)), n=1, \dots\}
$$

The map $\gz$ and the invariant measure $\nu$ are extended to ${\overline \Delta}$ in the natural way.
We shall still denote the probability with respect to the extended measure by $\Prob$.

We extend the map $\Phi$ to a map 
$${\overline \Phi}:{\overline \Delta}\to \Omega_{\A,B}^{\mathbb Z}, $$

$$
{\overline \Phi}({\bf \la})=\dots \omega_{-n}\dots \omega_0\dots \omega_n\dots,   
$$ 

if $(\la(-n), \pi(-n))\in\Delta(\omega_{-n})$, 
and $\gz^n(\la(0), \pi(0))\in\Delta(\omega_n).$

Now take a zippered rectangle $x\in\Omega(\R)$, $x=(\la, h, a, \pi)$. 
Set $\F^n(x)=(\la(n), h(n), a(n), \pi(n))$.

Consider a map 

\begin{equation}
\label{coding}
{\tilde \Phi}:\Y(\R)\rightarrow \Omega_{\A,B}^{{\mathbb Z}},
\end{equation}

given by 

$$
(\la,h,a,\pi)\ \rightarrow \dots \omega_{-n}\dots 
\omega_0\dots \omega_n\dots,
$$

where 

$$
(\la(n), \pi(n))\in\Delta(\omega_n)
$$

for all $n\in{\mathbb Z}$.

Under the natural projection $(\la,h,a,\pi)\rightarrow(\la,\pi)$, the $\F$-invariant measure 
${\overline \nu}$ on $\Y^{\pm}(\R)$ is mapped to the $\gz$-invariant measure $\nu$ on $\Delta(\R)$, whence 
the measure ${\tilde \Phi}_* {\overline \nu}$ is exactly the probability measure $\Prob$ on the space of bi-infinite 
sequences. To complete the identification of the spaces $(\Y^{\pm}(\R),{\overline \nu}) $ 
and $(\Omega_{\A,B}^{{\mathbb Z}}, \Prob)$, it remains to show that almost surely there is at most one
zippered rectangle corresponding to a given symbolic sequence.

\begin{proposition}
Let $\q\in \wab$ be such that all entries of the matrix $A(\q)$ are positive. 
Let $\omega\in  \Omega_{\A,B}^{{\mathbb Z}}$ be such that the word $\q$ occurs 
infinitely many times in $\omega$. Then there exists 
at most one zippered rectangle corresponding to $\omega$. 
\end{proposition}

Proof. Write 
$$
\omega=\dots \omega_{-n}\dots \omega_0\dots \omega_n\dots, 
$$

and let $(\la,h, a, \pi)$ be a zippered rectangle corresponding to $\omega$; we want to show 
that $(\la,h, a, \pi)$ is uniquely defiend by $\omega$.

First, $(\la,\pi)$ is uniquely defined by the "future" $\omega_0\dots \omega_n\dots$
of $\omega$.

Denote $w(n)=\omega_{-n}\dots \omega_0$, $(\la(-n), h(-n), a(-n), \pi(-n)=\F^{-n} (\la,h,a,\pi)$.

For any $n$ , the interval exchange $(\la(-n), \pi(-n))$ corresponds to the symbolic sequence 
$\omega_{-n}\dots \omega_0\dots$, and, again, is uniquely defined by that sequence.

By definition of the map $\F$, we have 

$$
\la(-n)=\frac{A(w(n))\la}{|A(w(n))\la|}, \ 
h(-n)=(A(w(n))^t)^{-1}h\cdot |A(w(n)\la|.
$$

Projectively, therefore, we have  

$$
{\mathbb R}_+h\subset A(w(n))^t{\mathbb R}_+^m. 
$$

Since the subword $\q$ occurs infinitely many times, the intersection

$$
\bigcap_{n=1}^{\infty} A(w(n))^t{\mathbb R}_+^m
$$

consists of a single line and the 
vector $h$ is therefore uniquely determined by 
the condition $<\la, h>=1$.

It remains to determine the vector $a$. 

By definition of the map $\F$, for any 
$n$ there exists an orthogonal matrix $U(-n)$, uniquely determined 
by $\omega$, and a vector $v(-n)$, 
uniquely determined by the the vectors $h(-n), \dots, h(0)$ and $\omega$, such that 

\begin{equation}
\label{formulafora}
\frac{U(-n)a(-n)+v(-n)}{|A(w(n)\la|}=a.
\end{equation}

Now let  $n$ be a moment such that all $\la(-n)_i>\frac 1{100m}$
(there are infinitely many such moments).  
Then $|a(-n)_i|<100m$ for all $i=1, \dots, m$ and,  (\ref{formulafora})
since $|A(w(n))\la|\to\infty$  as $n\to\infty$, 
(\ref{formulafora}) implies that $a$ 
is also uniquely determined by $\omega$.

The proof is complete.

\subsection{Transition probabilities.}

Take a sequence $c_1 \dots c_n \dots \in\Omega_{\A,B}$. 
Following Sinai \cite{sinailectures}, consider the {\it transition probability}

$$
\Prob(\omega_1=c_1|\omega_2=c_2, \dots, \omega_n=c_n, \dots)=\lim_{n\to\infty} \frac{\Prob (c_1c_2\dots c_n)}{\Prob (c_2\dots c_n)}. 
$$

In this subsection, we give a formula for this probability in terms of $(\la,\pi)=\Phi^{-1}(c_2\dots c_n\dots)$.

Assume $w_1\in\A$ is compatible with $(\la,\pi)$.

Denote  
$$
\Prob(w_1|(\la,\pi))=\Prob(((\la(-1), \pi(-1))=t_{w_1}(\la(0),\pi(0))|(\la(0),\pi(0))=(\la,\pi)).
$$

If $w_1\in\A$ is compatible with $(\la,\pi)$, from the definition of $\gz$ and from (\ref{detja}) we have 

\begin{equation}
\label{probwone}
\Prob(w_1|(\la,\pi))=\frac{\rho(t_{w_1}(\la,\pi))}{\rho(\la,\pi)|A(w_1)\la|^m}
\end{equation}

Since the invariant density is a homogeneous function of degree $-m$, we have 
$$
\rho(T_{w_1}(\la,\pi))=\frac{\rho(t_{w_1}(\la,\pi))}{|A(w_1)\la|^m},
$$
and we can rewrite (\ref{probwone}) as follows:

\begin{equation}
\label{probwoneunif}
\Prob(w_1|(\la,\pi))=\frac{\rho(T_{w_1}(\la,\pi))}{\rho(\la,\pi)}
\end{equation}

Let $w=w_1\dots w_n$ be compatible with $(\la,\pi)$.

Denote  
$$
\Prob(w|(\la,\pi))=\Prob((\la(-k), \pi(-k))=t_{w_{n-k+1}}(\la(1-k),\pi(1-k)), k=1, \dots,n|(\la(0), \pi(0))=(\la,\pi)).
$$

From (\ref{probwone}), by induction, we have 

\begin{equation}
\label{trpr}
\Prob(w|(\la,\pi))=\frac{\rho(t_w(\la,\pi))}{\rho(\la,\pi)|A(w)\la|^m}
\end{equation}

Since the invariant density is a homogeneous function of degree $-m$, we have 
$$
\rho(T_{w}(\la,\pi))=\frac{\rho(t_{w}(\la,\pi))}{|A(w)\la|^m},
$$
and we can rewrite (\ref{trpr}) as follows:

\begin{equation}
\label{homdens}
\Prob(w|(\la,\pi))=\frac{\rho(T_w(\la,\pi))}{\rho(\la,\pi)}
\end{equation}

\begin{corollary}
\label{bounddens}
There exists $C>0$ such that the following is true.
Suppose  $w\in W_{\A,B}$ is compatible with $(\la,\pi)$. Then 

$$
\Prob(w|(\la,\pi))\geq \frac C{\rho(\la,\pi)|A(w)|^m}
$$
\end{corollary}

Proof: recall that the invariant density is a positive homogeneous function of degree $-m$ 
and therefore  is bounded from below: there exists $C>0$ such that 
$\rho(\la,\pi)>C$ for all $(\la,\pi)\in\Delta(\R)$. In particular, 
$\rho(t_w(\la,\pi))>C$. Substituting into (\ref{trpr}), we obtain the result.

For $\epsilon:0<\epsilon<1$,  let 
$$
\Delta_{\epsilon}=\{(\la,\pi)\in\Delta(\R), \min|\la_i|\geq \epsilon\}.
$$

For any $\epsilon>0$ there exists a constant $C(\epsilon)$ such that for any  $(\la,\pi)\in\Delta_{\epsilon}$ we have 
$\rho(\la,\pi)<C(\epsilon)$. 
\begin{corollary} 
\label{epsilonbound}
For any $\epsilon>0$ there exists  $C(\epsilon)>0$ such that 
if $(\la,\pi)\in\Delta_{\epsilon}$, then 

$$
\Prob(w|(\la,\pi))\geq \frac {C(\epsilon)}{|A(w)|^m}.
$$
\end{corollary}

\section{Proof of the Exactness}
\label{proofexactness}

First, one notes that the discrete parameter $\pi$ does not give rise to any period, and then
the proof follows the standard pattern \cite{csf, viana}: since almost every point of any measurable subset is a density point, 
bounded distortion estimates of Proposition \ref{bdddist} imply that if the measure of a tail event is positive, then it
must be arbitrarily close to $1$.

In more detail, observe that there exists an integer $M$ such that for any $n>M$ and for any $\pi, \pip\in\R$ there 
exist $k_1, \dots, k_{2n}$ such that 
$a^{k_1}b^{k_2}\dots a^{k_{2n-1}}b^{k_{2n}}\pi=\pip$. This follows from conmnectedness of the Rauzy graph and 
the fact that for any $\pi\in\R$ there exist $n_1, n_2$ such that $a^{n_1}\pi=b^{n_2}\pi=\pi$.

Let $\alpha_0$ be the partition of $\Delta^+$ into $\Delta_{\pi}^+$, $\pi\in\R$, and let 
$\alpha_n$ be the partition into the cylinders $\Delta(w)$, where $w\in \W_{\A, B}$, $|w|=2n$.

\begin{lemma}
There exists $k>0$ such that the following is true.
Suppose $C\subset\Delta^+$, and there exists $\pi\in\R$ such that 
$\Delta_\pi^+\subset C$. Then $\gz^{2k}C=\Delta^+(\R)$. 
\end{lemma}

This implies 

\begin{lemma}
\label{contexact}
There exists $k>0$ such that the following holds. For any $\varepsilon>0$ there is $\delta>0$ such that 
for any  $C\subset \Delta^+(\R)$  satisfying ${\bf m}(C\bigtriangleup\Delta_{\pi}^+)<\delta$, we have
${\bf m}(\gz^{2k}C\bigtriangleup\Delta^+)<\varepsilon$. 
\end{lemma}

Now suppose $C\subset \Delta^+$ is a $\gz^2$-tail event, i.e., for any $n>0$ there exists 
$B_n$ such that $C=\gz^{-2n}B_n$ and $0<\nu(C)<1$. Then $\nu(B_n)=\nu(C)$ and, by 
Lemma \ref{contexact}, we can assume that there exists $\varepsilon>0$ such that
for any $\pi\in\R$, we have 

\begin{equation}
\label{epsmeas}
{\bf m}((\Delta^+\setminus C)\cap \Delta_{\pi}^+)\geq \varepsilon
\end{equation}

Let $\q=q_1\dots q_l$ be a word such that the matrix $A(\q)$ is positive.

For almost any $(\la,\pi)\in C$ we have

\begin{equation}
\label{densepoint}
\lim_{n\to\infty}\frac{{\bf m}(\alpha_n(\la,\pi)\cap C)}{{\bf m} (\alpha_n(\la,\pi))}=1
\end{equation}

Now let $n$ be such that  $\gz^{2n}(\la,\pi)\in \Delta(\q)$. Denote 
$(\lap, \pip)=\gz^{2n}(\la,\pi)$. Let $A$ be the corresponding renormalization matrix, 
that is, $\la=J_A\lap$. Then $A=A_1A(\q)$ for some (unimodular nonnegative integer) matrix $A_1$.
We have $\alpha_n(\la,\pi)=J_A(\Delta_{\pip}^+)$. By Proposition \ref{bdddist}, from 
(\ref{epsmeas}), we deduce that there exists $\varepsilon^{\prime}$, not depending on $n$ such that 

$$
\frac{{\bf m}(\alpha_n(\la,\pi)\cap (\Delta^+\setminus C))}{{\bf m} (\alpha_n(\la,\pi))}\geq \varepsilon^{\prime}.
$$

Since, by ergodicity, for almost any $(\la,\pi)$ we can find infinitely many $n$ such that 
$\gz^{2n}(\la,\pi)\in\Delta(\q)$, we arrive at a contradiction with (\ref{densepoint}), which gives the exactness of 
$\gz^2$.

\section{The Main Lemma}

We shall suppose from now on that the Rauzy class $\R$ is fixed and will often 
suppress it from notation.

For $\epsilon:0<\epsilon<1$,  define, in the same way as above,

$$
\Delta_{\epsilon}=\{(\la,\pi)\in\Delta(\R), \min|\la_i|\geq \epsilon\}.
$$

\begin{lemma}
\label{logsteps}
There exist positive constants $\gamma,K,p$ 
such that the following is true for any $\epsilon>0$. 
Suppose $(\la, \pi)\in\Delta_{\epsilon}$. Then 

$$
\Prob \{\exists n\leq K|\log\epsilon|, (\la(-n),\pi(-n))\in\Delta_{\gamma}|
(\la(1), \pi(1))=(\la,\pi))\}\geq p.
$$
\end{lemma}

From Corollary \ref{epsilonbound}, we obtain

\begin{corollary}
\label{qlogsteps}
Let $\q\in W_{\A,B}$, $\q=q_1\dots q_l$ be such that 
all entries of the matrix $A(\q)$ are positive.
Then there exist positive constants $K(\q),p(\q)$ 
such that the following is true for any $\epsilon>0$. 
Suppose $(\la, \pi)\in\Delta_{\epsilon}$. Then 

$$
\Prob \{\exists n\leq K(\q)|\log\epsilon|, (\la(-n),\pi(-n))\in\Delta(\q)|
(\la(1), \pi(1))=(\la,\pi))\}\geq p(\q).
$$
\end{corollary}

Informally, the proof of Lemma \ref{logsteps} proceeds by getting rid 
of small intervals.

For $\gamma>0$, $k\leq m$, denote 

$$
\Delta_{\gamma, k}=\{(\la,\pi): \exists i_1, \dots, i_k: \la_{i_1}, \dots, \la_{i_k}\geq \gamma\}.
$$

and 

$$
\Delta_{\gamma, k, \epsilon}=\{(\la,\pi): \la_i\geq \epsilon {\rm \ for \  all\ } i=1, \dots, m {\rm \ and\ }
\exists i_1, \dots, i_k: \la_{i_1}, \dots, \la_{i_k}\geq \gamma \}.
$$

Lemma \ref{logsteps} follows from

\begin{lemma}
\label{onemore}
There exist constants $L,K,p$, depending only on the Rauzy class, such that 
the following is true for any $\gamma, k, \epsilon$. 

Assume $(\la,\pi)\in\Delta_{\gamma, k, \epsilon}$. 

Then 
$$
\Prob \{\exists n\leq K|\log\epsilon|:  (\la(-n),\pi(-n))\in\Delta_{\gamma/L, k+1, \epsilon/L}|
(\la(1), \pi(1))=(\la,\pi))\}\geq p.
$$
\end{lemma}

Lemma \ref{onemore} is proved in the next four sections.

\section{An estimate on the number of Rauzy operations.}

Recall that, if $(\la,\pi)\in\Delta^+$, then the $\gz$-preimages of $(\la,\pi)$ are 
the exchanges  $t_{a^{-n}}(\la,\pi)$, $n=1, \dots$.
whereas if $(\la,\pi)\in\Delta^-$, then  the $\gz$-preimages of $(\la,\pi)$ are 
the exchanges  $t_{b^{-n}}(\la,\pi)$, $n=1, \dots$.

Denote 

$$
{\bf p}_n(\la,\pi)=
\begin{cases}
\Prob((\la(-1), \pi(-1))=t_{a^{-n}}(\la,\pi)|(\la(0), \pi(0))=(\la,\pi)), &\text{if $(\la,\pi)\in\Delta^+$ ;}\\
  \Prob((\la(-1), \pi(-1))=t_{b^{-n}}(\la,\pi)|(\la(0), \pi(0))=(\la,\pi)), &\text{if $(\la,\pi)\in\Delta^-$.}
\end{cases}
$$

For $\la\in {\mathbb R}^m_+$, set
$$
T_{a^{-1}}^{(\pi)}(\la)=A(a^{-1}\pi,a)\la, \ t_{a^{-1}}^{(\pi)}(\la)=\frac{A(a^{-1}\pi,a)\la}{|A(a^{-1}\pi,a)\la|},
$$

$$
T_{b^{-1}}^{(\pi)}(\la)=A(b^{-1}\pi,b)\la, \ t_{b^{-1}}^{(\pi)}(\la)=\frac{A(b^{-1}\pi,b)\la}{|A(b^{-1}\pi,b)\la|},
$$
and 

$$
T_{a^{-n}}^{(\pi)}(\la)=T_{a^{-1}}^{(a^{1-n}\pi)}\dots T_{a^{-1}}^{(\pi)}\la,\  
t_{a^{-n}}^{(\pi)}(\la)=t_{a^{-1}}^{(a^{1-n}\pi)}\dots t_{a^{-1}}^{(\pi)}\la,
$$
$$
T_{b^{-n}}^{(\pi)}(\la)=T_{b^{-1}}^{(b^{1-n}\pi)}\dots T_{b^{-1}}^{(\pi)}\la, \ 
t_{b^{-n}}^{(\pi)}(\la)=t_{b^{-1}}^{(b^{1-n}\pi)}\dots t_{b^{-1}}^{(\pi)}\la,
$$

so that we have

$$
t_{a^{-n}}(\la,\pi)=(t_{a^{-n}}^{(\pi)}\la, a^{-n}\pi), \ T_{a^{-n}}(\la,\pi)=(T_{a^{-n}}^{(\pi)}\la, a^{-n}\pi),
$$
$$
t_{b^{-n}}(\la,\pi)=(t_{b^{-n}}^{(\pi)}\la, b^{-n}\pi), \ T_{b^{-n}}(\la,\pi)=(T_{b^{-n}}^{(\pi)}\la, b^{-n}\pi).
$$

\begin{lemma}
\label{tailestimate}
If $(\la,\pi)\in\Delta^+$, then, for any $N\geq 1$, we have

$$
\sum_{n=N+1}^{\infty}{\bf p}_n(\la,\pi)=\frac{\rho_{a^{-N}\pi}^+(T_{a^{-N}}^{(\pi)}(\la))}{\rho_{\pi}^+(\la)}
$$

If $(\la,\pi)\in\Delta^-$, then, for any $N\geq 1$, we have

$$
\sum_{n=N+1}^{\infty}{\bf p}_n(\la,\pi)=\frac{\rho_{b^{-N}\pi}^-(T_{b^{-N}}^{(\pi)}(\la))}{\rho_{\pi}^-(\la)}
$$

\end{lemma}

Proof: We only consider the case $(\la,\pi)\in\Delta^+$.
In this case, the formula (\ref{probwoneunif}) can be written as
$$
{\bf p}_n(\la,\pi)=\frac{\rho^-_{a^{-n}\pi}(T_{a^{-n}}^{(\pi)}\la)}{\rho^+_{\pi}(\la)},
$$ 

whence we can write

\begin{equation}
\label{invce}
\rho^+_{\pi}(\la)=\sum_{n=1}^{\infty}\rho^-_{a^{-n}\pi}(T_{a^{-n}}^{(\pi)}\la).
\end{equation}

Note that this formula is true for {\it any} permutation $\pi$ and {\it any} $\la$ 
(i.e., even if $\la\notin \Delta_{\pi}^+$, the formula, being an identity between rational functions, 
still holds).

Since, for any $\la$, we have 
$$
T_{a^{-n-N}}^{(\pi)}\la=T_{a^{-n}}^{(a^{-N}\pi)}(T_{a^{-N}}^{(\pi)}\la),
$$

from (\ref{invce}) we obtain

$$
\rho^+_{a^{-N}\pi}(T_{a^{-N}}^{(\pi)}\la)=\sum_{n=1}^{\infty}\rho^-_{a^{-n-N}\pi}T_{a^{-n-N}}^{(\pi)}\la=
\rho^+_{\pi}(\la)(\sum_{n=N+1}^{\infty}{\bf p}_n(\la,\pi)),
$$

and the Lemma is proved.

\subsection{Bounded growth}

Let $(\la, \pi)\in\Delta(\R)$.

Define 
$$
(\la^{(n)}, \pi^{(n)})=\begin{cases}
t_{a^{-n}}(\la,\pi)), &\text{if $(\la,\pi)\in\Delta^+$ ;}\\
t_{b^{-n}}(\la,\pi),  &\text{if $(\la,\pi)\in\Delta^-$.}
\end{cases}
$$

$$
(\La^{(n)}, \pi^{(n)})=\begin{cases}
T_{a^{-n}}(\la,\pi)), &\text{if $(\la,\pi)\in\Delta^+$ ;}\\
T_{b^{-n}}(\la,\pi),  &\text{if $(\la,\pi)\in\Delta^-$.}
\end{cases}
$$
We have  
$$
\gz^{-1}(\la,\pi)=\{(\la^{(n)}, \pi^{(n)}), n=1, \dots\}.
$$  
and
$$
{\bf p}_n=\Prob ((\la(-1), \pi(-1))=(\la^{(n)}, \pi^{(n)})|((\la(0), \pi(0))=(\la,\pi)).
$$

For any $n\in{\mathbb N}$, there exists $i(n)\in\{1, \dots, m\}$ such that  

$$
|\La^{(n)}|-|\La^{(n-1)}|=\la_{i(n)}. 
$$

If $(\la(-1), \pi(-1))$ is a $\gz$-preimage of $(\la,\pi)$ and $(\la(-1), \pi(-1)=t_{c^{-n}}(\la,\pi)$, 
$c=a$ or $b$, then we define a vector $\La(-1)$ by the relation $(\La(-1), \pi(-1)=T_{c^{-n}}(\la,\pi)$
(in other words, $(\La(-1), \pi(-1))$ is the Zorich preimage without normalization).

\begin{lemma}
\label{bddgrowth}
There exists a constant $C(\R)$, depending on the Rauzy class only, such that  
for any $(\la,\pi)\in\Delta(\R)$ we have 
$$
\Prob(|\La(-1)|>K|(\la(0), \pi(0)=(\la,\pi))<\frac {C(\R)}{K-2}.
$$
\end{lemma}

For definiteness, assume $\la\in\Delta_{\pi}^{-}$ (the proof is completely 
identical in the other case). 
Then $\gz$-preimages of $(\la,\pi)$ are
$(\la^{(n)},\pi^{(n)})=t_{b^{-n}}(\la,\pi)$, $n=1,2,\dots$.

By construction \cite{zorich},  the invariant density $\rho_{\pi}^-$ has the form
$$
\rho_{\pi}^{-}(\la)=\sum_{i=1}^N\frac 1{l_{i1}(\la)l_{i2}(\la)\dots l_{i m}(\la)},
$$
where the functions $l_{ij}$ are linear:
$$
l_{ij}(\la)=a_{ij}^{(1)}\la_1+\dots+a_{ij}^{(m)}\la_m,
$$

and all $a_{ij}^{(r)}$ are nonnegative (in fact, $a_{ij}^{(r)}=0$ or $1$, 
but we do not need this fact here).

Let $l$ be the length of the $a$-cycle of 
$\pi$, that is, the smallest such number that $a^l\pi=\pi$. 

Since for any $k>0$ we have $a^{-kl}\pi=\pi$, 
from Lemma \ref{tailestimate} we obtain

$$
\sum_{n= kl+1}^{\infty}{\bf p}_n(\la,\pi)=
\frac{\rho_{\pi}^{-}(\La^{(kl)})}{\rho_{\pi}^{-}(\la)}.
$$ 

As noted above, for any $n>0$ there exists $\la_{i(n)}$ such that 

$$
|\La^{(n)}|-|\La^{(n-1)}|=\la_{i(n)},
$$
and, in fact, 
$$
\La^{(n)}=(\la_1, \dots, \la_{m-1}, \la_m+\la_{i(1)}+\dots+\la_{i(n)}).
$$

Since 
$$
\sum_{n= kl+1}^{\infty}{\bf p}_n(\la,\pi)\to 0 {\rm \ as \ } k\to\infty,
$$

for any $i=1, \dots, N$ there exists $j$ such that $a_{ij}^{(m)}>0$. 
Renumbering, if necessary, the linear forms $l_{ij}$, we may assume that 
$a_{i1}^{(m)}>0$ for any $i$. Denote $\epsilon=\min a_{i1}^{(m)}$ and $L=\max a_{i1}^{(r)}$.
For any $\la\in{\mathbb R}^m_+$ we have then  
$$
\epsilon \la_m \leq l_{i1}(\la)\leq L|\la|,
$$
whence 

\begin{equation}
\label{maintailest}
\frac{\rho_{\pi}^{-}(\La^{(kl)})}{\rho_{\pi}^{-}(\la)}\leq  \frac{L}{\epsilon(\la_m+\la_{i(1)}+\dots+\la_{i(kl)})}.
\end{equation}

Let $N$ be the smallest number such that $|\La(-N)|>K$ and let $s$ be the largest such  integer that $sl<N$. 
Then $|\La(-sl)|>K-1$ (because all $\la_{i(sl+1)}, \dots, \la_{i(N)}$ are all distinct) and  
$\la_m+\la_{i(1)}+...\dots +\la_{i(sl)}>K-2$ (because  $|\La(-sl)|=1+
\la_{i(1)}+...\dots +\la_{i(sl)}$.

 Therefore, by (\ref{maintailest}), we obtain
$$
\frac{\rho_{\pi}^{-}(\La^{(kl)})}{\rho_{\pi}^{-}(\la)}\leq \frac{L}
{\epsilon}\frac 1{K-2},
$$
and the Lemma is proved.

\begin{lemma}
\label{trivest}
Suppose $(\la,\pi)\in\Delta^+$, and let $l$ be the length of the $a$-cycle of $\pi$. 

Then, for any $k\geq 1$, we have 

$$
\sum_{n=kl+1}^{\infty}{\bf p}_n(\la,\pi)\geq (\frac{\la_{\pi^{-1}m}}{\la_{\pi^{-1}m}+k})^m.
$$

Suppose $(\la,\pi)\in\Delta^-$, and let $l$ be the length of the $b$-cycle of $\pi$. 

Then, for any $k\geq 1$, we have 

$$
\sum_{n=kl+1}^{\infty}{\bf p}_n(\la,\pi)\geq (\frac{\la_{m}}{\la_{m}+k})^m.
$$
\end{lemma}

Proof. Again, we only consider the case $(\la,\pi)\in\Delta^-$, as the proof of the other case is identical.
$$
\sum_{n=kl+1}^{\infty}{\bf p}_n(\la,\pi)=
\frac{\rho_{\pi}^-(\La^{(kl)}))}{\rho_{\pi}^-(\la)}.
$$

Set $\La^{(kl)}=(\La^{(kl)}_1, \dots, \La^{(kl)}_m).$ 

For $k=1$ we have  $\La^{(l)}_i=\la_i$ for $i<m$ and $\La^{(l)}_m=\la_m+\la_{i(1)}+ \dots+ \la_{i(l)}$, and for arbitrary $k$
by induction we obtain
$\La^{(kl)}_i=\la_i$ for $i<m$ and $\La^{(kl)}_m=\la_m+k(\la_{i(1)}+ \dots+ \la_{i(l)}).$

Note that $\la_{i(1)}+\dots  \la_{i(l)}\leq 1$ (since $i(1), \dots, i(l)$ are all distinct).

As in the proof of the previous Lemma, write

$$
\rho_{\pi}^{-}(\la)=\sum_{i=1}^N\frac 1{l_{i1}(\la)l_{i2}(\la)\dots l_{i m}(\la)},
$$

whence 

\begin{equation}
\label{compratio}
\frac{\rho_{\pi}^-(\La^{(kl)}))}{\rho_{\pi}^-(\la)}\geq \min_i 
\frac{l_{i1}(\la)l_{i2}(\la)\dots l_{i m}(\la)}{l_{i1}(\La^{(kl)})l_{i2}(\La^{(kl)})\dots l_{i m}(\La^{(kl)})}. 
\end{equation}

For any linear form $l(\la)=a_1\la_1+\dots+a_m\la_m$, $a_i\geq 0$, we have 

$$
\frac{l(\La^{(kl)})}{l(\la)}\geq \frac{\la_m}{\la_m+k(\la_{i(1)}+ \dots+ \la_{i(l)}}\geq\frac{\la_m}{\la_m+k},
$$

and the Lemma follows.

\section{An estimate on the probability of stopping.}

\begin{lemma}
\label{addbig}

For any $\gamma>0$, there exists $c(\gamma)>0$ such that
if $\la_{i(N)}>\gamma$, then 

$$
\frac{{\bf p}_N(\la,\pi)}{\sum_{n=N+1}^{\infty}{\bf p}_n(\la,\pi)}\geq c(\gamma)
$$
\end{lemma}

From Lemma \ref{trivest} we immediately have the following Corollary.

\begin{corollary}
\label{coraddbig}

For any $\gamma>0$, there exists $c(\gamma)>0$ such that
the following is true.

Assume $(\la,\pi)\in\Delta^+$, $\la_{i(N)}>\gamma$,  $\la_{\pi^{-1}m}>\gamma$. Then
$$
{\bf p}_N\geq \frac{c(\gamma)}{N^m}.
$$
Similarly, if $(\la,\pi)\in\Delta^-$, $\la_{i(N)}>\gamma$,  $\la_{m}>\gamma$, then  
$$
{\bf p}_N\geq \frac{c(\gamma)}{N^m}.
$$
\end{corollary}

If $(\la,\pi)\in\Delta^+$, then, by the definition of ${\bf p}_n(\la,\pi)$ and by Lemma \ref{tailestimate}, we have 
$$
{\bf p}_N(\la,\pi)=\frac{\rho^-_{a^{-N}\pi}(T_{a^{-N}}^{(\pi)}\la)}{\rho^+_{\pi}(\la)},
$$
$$
\sum_{n=N+1}^{\infty}{\bf p}_n(\la,\pi)=\frac{\rho_{a^{-N}\pi}^+(T_{a^{-N}}^{(\pi)}(\la))}{\rho_{\pi}^+(\la)},
$$
and, therefore, 
$$
\frac{{\bf p}_N(\la,\pi)}{\sum_{n=N+1}^{\infty}{\bf p}_n(\la,\pi)}=\frac{\rho^-_{a^{-N}\pi}(T_{a^{-N}}^{(\pi)}\la)}
{\rho_{a^{-N}\pi}^+(T_{a^{-N}}^{(\pi)}(\la))}.
$$

Lemma \ref{addbig} follows now from the following 

\begin{lemma}
\label{abest}
For any $\gamma>0$ there exists a constant $c(\gamma)>0$ such that the following is true.
Let $(\la,\pi)\in\Delta(\R)$. If $\la_{\pi^{-1}m+1}>\gamma$, then 
$$
\frac{\rho^-_{\pi}(\la)}{\rho^+_{\pi}(\la)}\geq c(\gamma).
$$
If $\la_{\pi^{-1}(\pi(m)+1)}>\gamma$, then 
$$
\frac{\rho^+_{\pi}(\la)}{\rho^-_{\pi}(\la)}\geq c(\gamma).
$$
\end{lemma}

The proof of Lemma \ref{abest} will take the remainder of this section.

First, we modify Veech's coordinates on the space of zippered rectangles.  
Take a zippered rectangle $(\la, h, a, \pi)\in\Delta(\R)$, and introduce the vector
$\delta=(\delta_1, \dots, \delta_m)\in {\mathbb R}^m$ by the formula

$$
\delta_i=a_{i-1}-a_i, \ i=1, \dots, m
$$

(here we assume, as always, $a_0=a_{m+1}=0$).

\begin{proposition}
\label{deltazip}
 The data $(\la,\pi,\delta)$ determine the zippered rectangle $(\la,h, a,\pi)$ uniquely.
\end{proposition}

{\bf Remark.} The coordinates  $(\la,\pi, \delta)$ on the space of zippered rectangles
have a natural interpretation in terms of the cohomological coordinates 
of Hubbard and Masur \cite{hubmas}: namely, the $\la_i$ are the real parts of the 
corresponding cycles, and the $\delta_i$ are (minus) the imaginary parts.

Proof of Proposition \ref{deltazip}. 
For any $i=1, \dots, m$, we have  

\begin{equation}
\label{adelta}
a_i=-\delta_1-\dots-\delta_i,
\end{equation}

so the vector $a$ is uniquely defined by $\delta$.
It remains to show that the vector $h$ is uniquely defined by $\delta$, and, to do this, 
we shall express the $h$ through the $a$. 
First note that 
$$
h_{\pi^{-1}m}=a_{\pi^{-1}m}-a_m.
$$

Now, if $i\neq \pi^{-1}m$, then $i=\pi^{-1}(k-1)$ for some 
$k\in\{1, \dots, m\}$. The equation
\begin{equation}
h_i-a_i=h_{\pi^{-1}(\pi(i)+1)}-a_{\pi^{-1}(\pi(i)+1)-1}.
\end{equation}
then takes the form 
$$
h_{\pi^{-1}(k-1)}-a_{\pi^{-1}(k-1)}=h_{\pi^{-1}(k)}-a_{\pi^{-1}(k)-1},
$$
or, equivalently, 
$$
h_{\pi^{-1}(k)}=a_{\pi^{-1}(k)-1}+h_{\pi^{-1}(k-1)}-a_{\pi^{-1}(k-1)}.
$$

Since 
$$
h_{\pi^{-1}1}=a_{\pi^{-1}1-1},
$$
by induction, we obtain
$$
h_{\pi^{-1}k}=a_{\pi^{-1}k-1}+\sum_{l=1}^{k-1}(a_{\pi^{-1}l-1}-a_{\pi^{-1}l})
$$

for any $k=1, \dots, m$, and the Lemma is proved. 

The above computations give us the following expression
for $h$ in terms of $\delta$:

\begin{equation}
h_{\pi^{-1}k}=
-\sum_{i=1}^{\pi^{-1}k-1} \delta_i+\sum_{l=1}^{k-1}\delta_{\pi^{-1}(l)}
\end{equation}

or, equivalently,

\begin{equation}
h_{r}=-\sum_{i=1}^{r-1} \delta_i+\sum_{l=1}^{\pi(r)-1}\delta_{\pi^{-1}l}.
\end{equation}

Rewriting the inequalities defining the zippered rectangle in terms 
of $\delta$, we obtain by a straightforward computation 
the following system: 
$$
\delta_1+\dots +\delta_i\leq 0,\ \  i=1, \dots, m-1.
$$

$$
\delta_{\pi^{-1}1}+\dots+\delta_{\pi^{-1}i}\geq 0, \ \  i=1, \dots, m-1.
$$

The parameter $a_m=-(\delta_1+\dots+\delta_m)$ can be both positive 
and negative.

Introduce the following cones in ${\mathbb R}^m$:

$$
K_{\pi}=\{\delta=(\delta_1,\dots,  \delta_m): \delta_1+\dots +\delta_i\leq 0,
\delta_{\pi^{-1}1}+\dots+\delta_{\pi^{-1}i}\geq 0, i=1, \dots, m-1\},
$$
$$
K_{\pi}^+=K_{\pi}\cap \{\delta:\sum_{i=1}^m\delta_i\leq 0\},
K_{\pi}^-=K_{\pi}\cap \{\delta:\sum_{i=1}^m\delta_i\geq 0\}.
$$

We have established the following 

\begin{proposition}
For $(\la,\pi)\in \Delta(\R)$ and an arbitrary $\delta\in K_{\pi}$
there exists a unique zippered rectangle $(\la,h,a,\pi)$ 
corresponding to the parameters $(\la,\pi,\delta)$.
\end{proposition}

In what follows, we shall simply refer to the zippered rectangle $(\la,\pi,\delta)$.

{\bf Remark.} It would be interesting to write down explicitly the genrating vectors for the cones
$K_{\pi}$, $K_{\pi}^+$, $K_{\pi}^-$; in particular, that would allow to give an explicit expression 
for the invariant densities of Veech \cite{veech} and Zorich \cite{zorich}.

Denote by $Area(\la,\pi,\delta)$ the area of the zippered 
rectangle $(\la,\pi, \delta)$. We have:

$$
Area(\la,\pi,\delta)=
\sum_{r=1}^m \la_rh_r=
\sum_{r=1}^m\la_r(-\sum_{i=1}^{r-1} \delta_i+
\sum_{l=1}^{\pi(r)-1}\delta_{\pi^{-1}l})=
$$
\begin{equation}
\label{ziparea}
\sum_{i=1}^m \delta_i (-\sum_{r=i+1}^m \la_r+\sum_{r=\pi(i)+1}^{m} 
\la_{\pi^{-1}r})=1.  
\end{equation}

A straightforward computation shows that
in the coordinates $(\la,\pi, \delta)$ the Rauzy induction map is written 
as follows: 

$$
{\cal T}(\la,\pi, \delta)=\begin{cases}
(\frac{A(\pi, b)^{-1}\la}{|A(\pi, b)^{-1}\la|}, b\pi, 
A(\pi, b)^{-1}\delta\cdot |A(\pi, b)^{-1}\la|), 
&\text{if $\la\in\Delta_{\pi}^+$;}\\
(\frac{A(\pi, a)^{-1}\la}{|A(\pi, a)^{-1}\la|}, a\pi, 
A(\pi, a)^{-1}\delta\cdot |A(\pi, a)^{-1}\la|), 
&\text{if $\la\in\Delta_{\pi}^-$.}
\end{cases}
$$

For $\la\in{\mathbb R}^m_+$, denote 

$$
K(\la,\pi)=K_{\pi}\cap \{\delta: Area(\la,\pi, \delta)\leq 1\},
$$
$$
K^+(\la,\pi)=K_{\pi}^+\cap \{\delta: Area(\la,\pi, \delta)\leq 1\},
$$
$$
K(\la,\pi)=K_{\pi}^-\cap \{\delta: Area(\la,\pi, \delta)\leq 1\}.
$$

Denote by $vol_m$ the Lebesgue measure in ${\mathbb R}^m$.

Set 
$$
{\bf r}(\la,\pi)=vol_m(K(\la,\pi)),
{\bf r}^+(\la,\pi)=vol_m(K^+(\la,\pi)),
{\bf r}^-(\la,\pi)=vol_m(K^-(\la,\pi)).
$$

By definition, the functions ${\bf r}, {\bf r}^+, {\bf r}^-$
are positive rational functions, homogeneous  of degree $-m$. 

\begin{lemma}
\begin{enumerate}
\item ${\bf r}^-(\la,\pi)={\bf r}(T_{b^{-1}}(\la,\pi))$.
\item ${\bf r}^+(\la,\pi)={\bf r}(T_{a^{-1}}(\la,\pi))$.
\item ${\bf r}(\la,\pi)={\bf r}(T_{a^{-1}}(\la,\pi))+
{\bf r}(T_{b^{-1}}(\la,\pi))$.
\end{enumerate}
\end{lemma}

Proof.  If 
$$
\delta=(\delta_1, \dots, \delta_m)\in K^-(\la,\pi),
$$ 
then
$$
{\tilde \delta}=(\delta_1, \dots, \delta_{m-1}, \delta_m+\delta_{\pi^{-1}m})\in K(T_{b^{-1}}(\la,\pi)),
$$ 
and vice versa. This gives a volume-preserving bijection between $K^-(\la,\pi)$ 
and $K(T_{b^{-1}}(\la,\pi))$, whence 
${\bf r}^-(\la,\pi)={\bf r}(T_{b^{-1}}(\la,\pi))$. 
The second assertion is proved in the same way, and the third 
follows from the first two.

\begin{corollary}
\label{rinvdensity}
$$
{\bf r}^+(\la,\pi)=\sum_{n=1}^{\infty} {\bf r}^-(T_{a^{-n}}(\la,\pi)).
$$
$$
{\bf r}^-(\la,\pi)=\sum_{n=1}^{\infty} {\bf r}^+(T_{b^{-n}}(\la,\pi)).
$$
\end{corollary}

We only prove the first assertion. We have 
$$
{\bf r}^+(\la,\pi)={\bf r}(T_{a^{-1}}(\la,\pi))={\bf r}^+(T_{a^{-1}}(\la,\pi)+{\bf r}^-(T_{a^{-1}}(\la,\pi)=
{\bf r}(T_{a^{-2}}(\la,\pi))+{\bf r}^-(T_{a^{-1}}(\la,\pi)).
$$

Proceeding by induction, 

$$
{\bf r}^+(\la,\pi)=\sum_{n=1}^{N} {\bf r}^-(T_{a^{-n}}(\la,\pi))+{\bf r}(T_{a^{-N-1}}(\la,\pi)).
$$

Since $$
T_{a^{-N-1}}(\la,\pi)=(T_{a^{-N-1}}^{(\pi)}(\la), a^{-N-1}\pi),
$$ 
and $|T_{a^{-N-1}}^{(\pi)}(\la)|\to\infty$ as $N\to\infty$, we obtain
${\bf r}(T_{a^{-N-1}}(\la,\pi))\to 0$ as $N\to\infty$, and the Corollary is proved.

Since the functions ${\bf r}, {\bf r}^+, {\bf r}^-$
are positive, rational and homogeneous  of degree $-m$,  Corollary \ref{rinvdensity}
implies that, for some positive constant $C(\R)$, depending only on the Rauzy class $\R$, we have 
$$
\rho^+(\la,\pi)=C(\pi){\bf r}^+(\la,\pi),
\rho^-(\la,\pi)=C(\pi){\bf r}^-(\la,\pi).
$$ 

By construction, for any $\la\in {\mathbb R}^m_+$ we have 
$$
{\bf r}_{\pi}^+(\la_1, \dots \la_m)={\bf r}_{\pi^{-1}}^-(\la_{\pi(1)}, \dots \la_{\pi(m)}).
$$

In view of this observation, it suffices to prove only the first assertion 
of the Lemma \ref{abest}, as the second one follows automatically.

Take $\delta=(\delta_1, \dots, \delta_m)\in {\mathbb R}^m$, and, for 
$\theta>0$, define 

$$
J^{(m)}_{\theta}\delta=(\delta_1, \dots, \delta_m+\theta),\ 
J^{(\pi^{-1}m)}_{-\theta}\delta=(\delta_1, \dots, 
\delta_{\pi^{-1}m}-\theta, \dots, \delta_m).
$$

\begin{proposition}
Let $\theta>0$. 
If $\delta\in K_{\pi}$, then $\jmtd\in K_{\pi}, \jmtpd\in K_{\pi}$.
If $\delta\in K_{\pi}^-$, then $\jmtd\in K_{\pi}^-$.
If $\delta\in K_{\pi}^+$, then $\jmtpd\in K_{\pi}^+$.
\end{proposition}
This follows directly from the definition of the cones $K_{\pi}, K_{\pi}^-, K_{\pi}^+$. 
From (\ref{ziparea}) we obtain 
 
$$
Area(\la,\pi, \jmtd)=Area(\la,\pi,\delta)+\theta(\sum_{r=\pi(m)+1}^m \la_{\pi^{-1}r}),
$$
$$
Area(\la,\pi, \jmtpd)=Area(\la,\pi,\delta)+\theta(\sum_{r=\pi^{-1}(m)+1}^m \la_r),
$$

which implies

\begin{proposition}
\label{comparea}
$$
Area(\la,\pi,\delta)\leq Area(\la,\pi, \jmtd)\leq Area(\la,\pi,\delta)+
\theta|\la|.
$$
$$
Area(\la,\pi,\delta)\leq Area(\la,\pi, \jmtpd)\leq Area(\la,\pi,\delta)+
\theta|\la|.
$$
\end{proposition}

For $s\in{\mathbb R}$ and a hyperplane of the form 
$\delta+\dots+\delta_m=s$, let $vol_{m-1}$ stand for the induced 
$(m-1)$-dimensional volume form on the hyperplane.

Denote 
$$
K_{s,\pi}=K_{\pi}\cap \{\delta:\sum_{i=1}^m\delta_i=s\},
$$
$$
K_s(\la,\pi)=K(\la,\pi)\cap K_{s,\pi},
$$
$$
V_s(\la,\pi)=vol_{m-1}(K_s(\la,\pi)).
$$

Denote by $\aminus$ the maximal possible value of $\delta_1+\dots+\delta_m=-a_m$ 
in $K(\la,\pi)$.

\begin{proposition}
\label{fromabove}
Assume $0\leq s\leq\aminus$. Then 
$$
V_s(\la,\pi)\leq (1+s)^{m-1}V_0(\la,\pi).
$$
\end{proposition}

Proof: Indeed, if $
(\la,\pi, \delta)\in V_s(\la,\pi),$ 
then Proposition \ref{comparea} implies   
$$(\la,\pi, \frac{J^{(\pi^{-1}m)}_{-s}\delta}{1+s})\in V_0(\la,\pi),$$ 
and the assertion follows.

\begin{proposition}
\label{frombelow}

Assume $s$, $0\leq s\leq 1$  is such that $\frac{s}{1-s}\leq \aminus$. 
Then 

$$
(\frac1{1-s})^{m-1}V_s(\la,\pi)\geq V_0(\la,\pi).
$$
\end{proposition}

Denote $\theta=\frac{s}{1-s}$, then $s=\frac{\theta}{1+\theta}$. 
If $(\la,\pi,\delta)\in V_0(\la,\pi),$ then 
$$(\la,\pi, \frac{\jmtd}{1+\theta})\in V_{s}(\la,\pi),
$$ and, again, the assertion follows.

Propositions \ref{fromabove}, \ref{frombelow} imply 

\begin{lemma}

For any $C_1>0$ there exists $C_2>0$ such that the following is true.

Let $\aminus(\la,\pi)<C_1$. Then 

$$
{\bf r}^-(\la,\pi)<C_2V_{m-1}^0(\la,\pi).
$$
\end{lemma}

Note that there exists $\epsilon>0$, depending only on $\R$ and  
such that for any $(\la,\pi)\in\Delta(\R)$, we have 
$\aminus>\epsilon$. In conjunction with Propositions \ref{fromabove}, 
\ref{frombelow}, this implies 
\begin{lemma}
There exists a constant $C_3$ such that for 
any $(\la,\pi)\in\Delta({\cal R})$, we have 

$$
\rho^-(\la,\pi)\geq C_3 V_0(\la,\pi).
$$
\end{lemma}

Since $a_m\leq h_{\pi^{-1}m+1}$, we have 

$$
\aminus(\la,\pi)\leq \frac 1 {\la_{\pi^{-1}m+1}}, 
$$
 which implies the following 
\begin{corollary}
For any $C_4>0$ there exists $C_5>0$ such that the following is true.

Assume $\la_{\pi^{-1}m+1}>C_4$. 
Then 
$$
\frac{{\bf r}^-(\la,\pi)}{{\bf r}^+(\la,\pi)}<C_5,
$$
\end{corollary}

which implies Lemma \ref{abest}.

\section{Kerckhoff names}

In the following two sections, we shall use Kerckhoff's  
convention of numbering the subintervals of an interval exchange \cite{kerckhoff}; 
to avoid confusion, we shall speak of {\it Kerckhoff names} of subintervals.

Take an interval exchange $(\la,\pi)$. A {\it Kerckhoff naming}  
on the subintervals
of $(\la,\pi)$ is defined by an arbitrary  permutation ${i_1, \dots, i_m}$ of the symbols $\{1, \dots, m\}$.
Once such a permutation is given,  we asign names $I_{i_1}, \dots, I_{i_m}$ 
to the subintervals of $(\la,\pi)$, from the left to the right
(i.e., the subinterval $[0, \la_1)$ is named $I_{i_1}$, the subinterval $[\la_1, \la_1+\la_2)$ is named 
$I_{i_2}$ and so forth).

A Kerckhoff naming of the subintervals of $(\la,\pi)$  
induces a naming on the subintervals of ${\cal \T}(\la,\pi)$ in the following way.  
Assume $\la_m<\la_{\pi^{-1}m}$ and the Rauzy operation $a$ was applied to $(\la,\pi)$ in order to obtain 
${\cal \T}(\la,\pi)$. Then the subintervals of  ${\cal \T}(\la,\pi)$ are named, from the left to the right, 
by $I_{i_1}, \dots, I_{i_{\pi^{-1}m}}, I_{m}$, $I_{\pi^{-1}m+1}, \dots, I_{m-1}$. If 
$\la_m>\la_{\pi^{-1}m}$  and the Rauzy operation $b$ was applied, then 
the subintervals of  ${\cal \T}(\la,\pi)$ are just named, as before, by $I_{i_1}, \dots, I_{i_m}$, 
from the left to the right. Proceeding inductively, we obtain a naming for any $\gz^n(\la,\pi)$. 
Conversely, if we have a Kerckhoff naming of subintervals  of $(\la,\pi)$, then, for any 
word $w\in\wab$ compatible with $(\la,\pi)$, we automatically obtain a  
Kerckhoff naming on the subintervalsof $t_w(\la,\pi)$ and $T_w(\la,\pi)$.

Let $(\la,\pi)$ be an interval exchange with a Kerckhoff naming $I_{i_1}, \dots, I_{i_m}$. 
If $(\la,\pi)\in\Delta^+$, then we say that $I_{i_{\pi^{-1}m}}$ is the subinterval {\it in the 
critical position} (we shall also sometimes say ``in the $a$-critical position'').  
If $(\la,\pi)\in \Delta^-$, then we say that $I_{i_{m}}$ is the subinterval {\it in the 
critical position} (we shall also sometimes say ``in the $b$-critical position'').

\section{Exponential growth.}

Let ${\bf x}\in{\overline \Delta}$, that is, ${\bf x}=(\dots, (\la(-n), \pi(-n), \dots, (\la,\pi))$,
where, as usual, $\gz(\la(-n), \pi(-n))=(\la(1-n), \pi(1-n))$.
Define the words $w(n)$ by the relation $(\la(-n), \pi(-n))=t_{w(n)}(\la,\pi)$. Set  
$(\La(-n), \pi(-n))=T_{w(n)}(\la,\pi)$.

 \begin{lemma}
\label{expgrowth}
There exists $N$ such that the following is true. 
For any  ${\bf x}\in{\overline \Delta}(\R)$,
there exist $i_1, i_2\in\{1,\dots,m\}$ such that 

$$
\La(-N)_{i_1}+\La(-N)_{i_2}\geq 2(\la(0)_{i_1}+\la(0)_{i_2})
$$
\end{lemma}

Proof: 

Take a point $x\in{\overline \Delta}$, 

$$
x=(\dots, (\la(-n), \pi(-n)), \dots, (\la,\pi)).
$$

Give Kerckhoff names $I_1, \dots, I_m$ to the subintervals of the 
exchange $(\la,\pi)$ from the left to the right, so that the length of 
$I_i$ is $\la_i$. We thus automatically obtain a Kerckhoff naming for the 
subintervals of $((\la(-n), \pi(-n))$ for any $n$. 

Let $I_{j_n}$ be the critical subinterval for $(\la(-n), \pi(-n))$. 

Consider the infinite sequence 
\begin{equation}
\label{cycseq}
I_{j_1}\dots I_{j_n}\dots. 
\end{equation}
Note that $j_n\neq j_{n+1}$. 
A subword $I_{j_k}\dots I_{j_{k+l}}$ will be called a {\it simple cycle} 
if $I_{j_k}=I_{j_{k+l}}$ whereas $I_{j_k}, \dots I_{j_{k+l-1}}$ are 
all distinct. Naturally, $l\leq m$. 
There are finitely many possible simple cycles, therefore 
there exists $N$, depending only on $m$, such that for any 
word of length $N$ in the alphabet $\{I_1, \dots, I_m\}$, some 
simple cycle occurs at least $m$ times. Now take the word 

\begin{equation}
I_{j_1}\dots I_{j_N},
\end{equation}
the beginning of the sequence (\ref{cycseq}), and take a simple cycle 
which occurs  $m$ times, say 

\begin{equation}
\label{simcyc}
I_{l_1}\dots I_{l_r},
\end{equation}

Here, of course, $r\leq m$.
Now estimate the non-renormalized length of 
the subintervals $I_{l_1}, \dots, I_{l_r}$ ($r\leq m)$. In the beginning, these 
are $\la_{l_1}, \dots, \la_{l_r}$. The key observation is, as usual, 
that the interval in critical position at a 
given inverse Zorich step was, at the previous step, added to the 
previous critical interval. 
After the first occurrence of the cycle (\ref{simcyc}), therefore, the 
(non-normalized) length of $I_{l_1}$ is at least $\la_{l_1}+\la_{l_2}$, 
that of $I_{l_2}$ is at least $\la_{l_2}+\la_{l_3}$ and so forth. 
After the second occurrence of (\ref{simcyc}), the length of 
$I_{l_1}$ is at least $\la_{l_1}+\la_{l_2}+\la_{l_3}$, 
that of $I_{l_2}$ is at least
$\la_{l_2}+\la_{l_3}+\la_{l_4}$, and so forth. Finally, after the 
$r$-th occurrence of 
(\ref{simcyc}),  the length of 
$I_{l_1}$ is not less than 
$\la_{l_1}+\la_{l_2}+\dots +\la_{l_r}$, that is, not less than 
$2\la_{l_1}$, since $\la_{l_1}=\la_{l_r}$. The Lemma is proven.

\section{Proof of the Lemma \ref{onemore}}

An informal sketch of the proof of Lemma \ref{onemore}. 
One divides the subintervals into ``big'' ones and ``small'' ones: the aim is to obtain one more ``big'' interval. 
For this, one must first put a small subinterval into critical 
position. This is achieved by Lemma \ref{posprob}. 
In the previous ection, we have seen that the total length of the 
(non-renormalized) interval grows 
exponentially with the number of Zorich steps (with an exponent depending on $\epsilon$).
When the total length of the interval doubles, we obtain a new ``big'' subinterval. 

\subsection{Putting a small interval into critical position}

Take an interval exchange $(\la,\pi)$ and name the subintervals
$I_1, \dots, I_m$, from the
right to the left.

\begin{proposition}
Any interval can be put both in the $a$-critical
and in the $b$-critical position.
\end{proposition}

Proof: First note that if an interval can be put in the $a$-critical
position, then it can also be put into the $b$-critical position just by performing the entire $a$-cycle of 
the corresponding permutation.
Since the permutation is irreducible, it suffices to prove that, if
$I_i$ can be put into critical position, then also all $I_j$ for $j>i$.
To prove this, take the shortest word $w$ that puts $I_i$ into the
$a$-critical position. Then, in the preimage, all $I_j$,
$j>i$, still stand to the right of $I_i$,though perhaps in a different
order (because an inversion of order between
$I_i$ and $I_j$ can only happen once $I_i$ reaches the critical position).
Therefore, we can immediately place any of the $I_j$, $j>i$, into the
$b$-critical position, but then also into the $a$-critical position.

More precisely, pick a positive integer 
$k\leq m$ and a real $\gamma>0$. We say that we have a 
$(k, \gamma)$-{\it big-small decomposition} if the intervals of the exchange are divided 
into two groups: 
$I_{i_1}, \dots, I_{i_k}$, each of length at least $\gamma$, and the remaining ones 
(nothing is said about the length of the remaining ones). 

Under the Kerckhoff convention, a big-small 
decomposition of $(\la,\pi)$ is inherited by all $t_w(\la,\pi)$ (one just takes the intervals with the same names).

\begin{lemma}
\label{posprob}For any $\gamma>0$, 
there exist constants $p(\gamma), L(\gamma)$ such that the following is true.
Let $(\la,\pi)\in\Delta_{k,\gamma}$ with a fixed big-small decomposition.
Then there exists $w\in {\cal W}_{{\cal A}, B}$ such that 
\begin{enumerate}
\item $\Prob(w|\la,\pi)\geq p(\gamma)$.
\item  $|T_w(\la,\pi)|<L(\gamma)$.
\item the exchange $t_w(\la,\pi)$ has a small interval in critical position.
\end{enumerate}
 \end{lemma}

Proof: Take the shortest 
word (in terms of the number of Zorich operations) 
that puts a small interval into critical position. Among all such words, pick the one 
that involves the smallest number of Rauzy operations.
The length of this word, as well 
as the number of Rauzy operations involved, only 
depends on the Rauzy class.
At each intermediate Rauzy step, all 
subintervals following the critical one either in the preimage or in the image 
must be big, otherwise there would exist a shorter word placing a small interval 
into critical position. Therefore, by Lemma \ref{addbig} and the Corollary \ref{coraddbig}, 
the probability of each Zorich operation involved is bounded from below by a constant that only 
depends on $\gamma$. The Lemma is proved.

\subsection{Completion of the proof.}

Proof: Take any $x\in{\overline \Delta}$. Take the first $n$ such that 
$|\La(-n)|>2$. By Lemma \ref{expgrowth}, $n<K|\log\epsilon|$.
By Lemma \ref{bddgrowth}, with positive probability depending only on $M$, we can 
assume $|\La(-n)|<2M$.
Consider two cases:
\begin{enumerate}
\item  at all steps from $1$ to $n$, only small intervals were added between themselves.
\item at some step a large interval was added to a small one.
\end{enumerate}

Note, that since we start with a small interval in critical position, either one 
or the other case holds (for, in order that a small interval be added to a big interval, 
a big interval must first be placed into critical position, and for that it must first be added 
to a small one). 

In the first case, the lengths of all large intervals remain the same, and after 
renormalization at step $n$, each large interval has length at least $\gamma/2M$.
However, since $|\La(-n)|>2$, there must be another interval of length at least
$1/2mM$, and the Lemma is proved.

In the second case, let $n_1$ be the first moment, at which a big interval 
is added to a small one. Then $|\La(-n_1)|<2$, and, since at previous moments only small 
intervals were added between themselves, 
we have $k+1$ intervals of length at least $\gamma/2$, and the Lemma is proved completely.

\section{Return times for the Teichm{\"u}ller flow.}

We have in fact proven a stronger statement, namely, the following Lemma.

\begin{lemma}
\label{flowtime}
For any word $\q\in\wab$ such that all entries of the matrix $A(\q)$
are positive, there exist constants $K_0(\q) ,p(\q)$, depending only on
$\q$ and  such that the following is true.
For any $K\geq K_0$ and any $(\la,\pi)\in\Delta(\R)$,

$$
\Prob(\exists n: (\la(-n),\pi(-n))\in\Delta_{\q},
|\La(-n)|<K)|(\la,\pi))\geq p(\q)
$$
\end{lemma}

This statement has the following Corollary for the Teichm{\"u}ller flow
on the space of zippered rectangles.

Take an arbitrary word $\q=q_1\dots q_{2l+1}\in\wab$ such that all
entries of the matrix $A(q_1\dots q_l)$ are positive and all entries of the matrix $A(\q)$ 
are positive. As usually, set 
$$
\Delta_{\q}=\{(\la,\pi): \Phi(\la,\pi)=\omega_1\dots
\omega_n\dots, \omega_1=q_1, \dots ,\omega_{2l+1}=q_{2l+1}\}.
$$

Consider also the cylinder
$$
{\overline \Delta}_{\q}=
\{\omega\in\Omega_{\A,B}^{\mathbb Z}, \omega_{-l}=q_1, \dots,
\omega_l=q_{2l+1}\}.
$$

Consider the flow $P^t$ as a special flow over
${\overline \Delta}_{\q}$. Denote the roof function
of the flow by $\tau_{\q}$.

We shall now see that Lemma \ref{flowtime} implies 
\begin{corollary}
\label{intcoc}
There exists $\epsilon>0$ such that
$$
\int_{{\overline \Delta}_{\q}} \exp(\epsilon 
\tau_{\q}(\omega))d\Prob(\omega)<+\infty.
$$
\end{corollary}

Take $\omega\in \Omega_{\A,B}^{{\mathbb Z}}$.
As usually, set
$$
(\la(-n), \pi(-n))=\Phi^{-1}(\omega_{-n}\dots \omega_{0}\omega_1\dots),
(\La(-n), \pi(-n))=T_{\omega_{-n}\dots \omega_{-1}\omega_0}(\la(0), \pi(0)).
$$
Set $n_{\q}(\omega)$ to be the smallest $n$ such that
$$
\omega_{-n}=q_1, \dots, \omega_{-n+2l}=q_{2l+1}.
$$
Finally, set $L_{\q}(\omega)=\log|\La(-n_{\q}(\omega))|$.
Informally, $L_{\q}(\omega)$ is the ``Teichm{\"u}ller flow time"
it takes $\omega$ to reach ${\overline \Delta}_{\q}$.

To establish the Corollary \ref{intcoc}, it suffices to prove

\begin{proposition}
\label{expflow}
There exists $\epsilon>0$ such that
$$
\int_{\Omega_{\A,B}^{{\mathbb Z}}} 
\exp(\epsilon L_{\q}(\omega))d\Prob(\omega)<+\infty.
$$
\end{proposition}

Proof of Proposition \ref{expflow}.
Our main tool will be Lemma \ref{flowtime}.
Take a $K>K_0$ such that $1-p(\q)+\frac 1 K<1$.
Define a random time $k_1(\omega)$ to be
the first moment $n$ such that $|\La(-n)(\omega)|>K$.
Note that the map
$$
{\tilde \sigma}(\omega)\to\sigma^{-k_1(\omega)}(\omega)
$$
is invertible (here, as always, $\sigma$ is the shift on $\Omega_{\A,B}$).

Introduce a function $\eta: \Omega_{\A, B}^{{\mathbb Z}}\to{\mathbb N}$ by
the formula
$$
\eta(\omega)=[\frac{\log |\La(-k_1(\omega))|}{\log K}].
$$

In other words, $\eta(\omega)=n$ if
$$
K^n\leq |\La(-k_1(\omega))|\leq K^{n+1}.
$$

\begin{proposition}
There exists a constant $C$ such that the following is true for any 
$K>K_0$.

For any $c_1\dots c_n \dots \in \Omega_{\A,B}^+$,

$$
\Prob(\{\omega: \eta(\omega)=n|\omega_1=c_1, \dots \omega_n=c_n\dots )\}
\leq \frac C {K^n}.
$$
\end{proposition}

This immediately follows from Lemma \ref{bddgrowth}.

\begin{proposition}
$$
\Prob(\{\omega: \eta(\omega)=1,
\omega_{-k_1(\omega)}\dots \omega_0 \ {\rm does \ 
not \ contain \ the \ word} 
\ \q\})\leq 1-p(\q)
$$
\end{proposition}

Finally, take a large $N$ and let
$$
n_N(\omega)=\min n: k_1(\omega)+\dots k_1({\tilde \sigma}^{-n}(\omega))\leq 
N.
$$

Note that, by definition,
$$
K^N\leq |\La(-n_1)(\omega)|\leq K^{2N}.
$$
Now consider the set
$$
{\tilde \Omega(N)}=
\{\omega: \omega_{-n_N(\omega)}\dots \omega_{0}\  {\rm does \ not \
contain \ the
\ word \ } \q\}.
$$
Note that
$$
\{\omega: L_{\q}(\omega)>2N \}\subset {\tilde \Omega(N)}.
$$
It suffices, therefore, to prove that there exists $r<1$ such that
$$
\Prob({\tilde \Omega}(N))\leq r^N.
$$
But by the previous two propositions, we immediately have
$$
\Prob({\tilde \Omega})\leq C(1-p(\q)+\frac 1{K})^N,
$$
and, since $1-p(\q)+\frac 1{K}<1$, the Proposition follows.

This Proposition admits an equivalent formulation in
terms of the norms of renormalization matrices
on the space of of interval exchange transformations.

More precisely, for $(\la,\pi)\in\Delta_{\q}$, $\Phi(\la,\pi)=
\omega_1\dots \omega_n\dots$,
we let $n^{\q}(\la,\pi)$ to be the
smallest $n>0$ such that
$\gz^{n}(\la,\pi)\in\Delta_{\q}$, and we set

$$
{\cal N}(\la,\pi)=||A(\omega_1\dots \omega_{n^{\q}(\omega)})||.
$$
\begin{corollary}
\label{intcocycle}
There exists $\epsilon>0$ such that
$$
\int_{\Delta_{\q}} N(\la,\pi)^{\epsilon}d\Prob<+\infty.
$$
\end{corollary}

{\bf Remark.}
First results on exponential decay for the probabilities
of return times were obtained by Jayadev Athreya.
In his approach, Athreya used the dynamics of $SL(2, {\mathbb R})$-action,
which allowed him to obtain optimal exponents.
The argument above is an attempt to recover some of Athreya's 
theorems using the language of interval exchange transformations; 
the argument above does not, however, give an optimal exponent.

Avila, Gou{\"e}zel, and Yoccoz have recently announced exponential
decay of correlations for the Teichm{\"u}ller flow.
One of the steps in their proof is, again,
an exponential estimate for return times,
which they have obtained independently (Avila [oral communication]).
Their exponent is optimal.

\section{Estimate of the measure.}

\begin{lemma}

There exists a constant $C(\R)$ depending only on the Rauzy class $\R$ such that 

$$
\nu(\Delta(\R)\setminus \Delta_{\epsilon}(\R))<C\epsilon
$$
\end{lemma}

The proof repeats that of Proposition 13.2 in Veech \cite{veech}.

Lemma \ref{logsteps} and
Corollary \ref{qlogsteps} therefore imply the following 

\begin{corollary}
\label{sqroot}
Let $\q\in W_{\A,B}$, $\q=q_1\dots q_l$ be such that 
all entries of the matrix $A(\q)$ are positive.

There exist $C>0, \alpha>0$ such that the following is true for any $n$.

$$
\Prob((\la,\pi):\gz^{2k}(\la,\pi)\notin \Delta(\q) \ {\rm for \ all \ }\ k, 1\leq k\leq n)\leq C\exp(-\alpha \sqrt n).
$$

\end{corollary}

Proof: Let $n=r^2$ and 
denote 
$$
X(n,\q)=\{(\la,\pi):\gz^{2k}(\la,\pi)\notin \Delta(\q) \ {\rm for \ all \ }\ k, 1\leq k\leq n)\}.
$$

Take 
$$
B(n)=\{(\la,\pi): \gz^{2k}(\la,\pi)\notin \Delta_{\exp(-r)} \ {\rm for \ some \ }\ k, 1\leq k\leq n)\}
$$

Then, by the previous Lemma, $\nu(B(n))\leq Cr^2\exp(-r)$, whereas, by Corollary \ref{qlogsteps}, 
$$
\nu(X(n, \q)\setminus B(n))\leq (1-p(\q))^r,
$$

and Corollary \ref{sqroot} is proven.

{\bf Remark.} This result allows to use the tower method of L.-S. Young \cite{lsy} and 
to obtain the decay rate $\exp(-\alpha \sqrt n)$ 
for correlations of bounded H{\"o}lder functions. 
For bounded Lipschitz functions, one can also use the method of 
V. Maume-Deschamps \cite{maume} and obtain the uniform rate of decay 
at the rate $\exp(-\alpha n^{1/2-\epsilon})$. It is not clear to me, however, how to 
use either of these methods in the invertible case.

\section{Inequalities}

Let 

$$
W^+_{{\cal A}, B}=\{w\in W_{\A,B}:\  |w| \  {\rm is  \ even\ },\  \Delta(w)\subset \Delta^+\}.
$$

\begin{lemma}
For any $C_1, C_2>0$ 
there exists $C_3>0$ such that the following is true.

Suppose $row(A)<C_1$ and $\la\in\Delta_{C_2}$.

Then 
$$
\frac1{C_3}\leq \frac{|A\la|^m}{\Pi_{j=1}^m\sum_{i=1}^mA_{ij}}\leq C_3
$$

\end{lemma}

Proof: 

Denote $A_j=\sum_{i=1}^mA_{ij}$, so that $|A|=\sum_{j=1}^m A_j$.

Then 
$$
\frac{A_j}{A_k}\leq row(A),
$$
whence 
$$
\frac {A_j}{|A|}\geq \frac 1 {m \  row(A)}.
$$

Finally, if $\la\in\Delta_{C_2}$, then 
$$
|A\la|\geq C_2|A|,
$$

which completes the proof.

\begin{corollary}
For any $C_4>0$, $C_5>0$ there exists $C_6>0$ such  that the following is true.
Suppose $(\la,\pi)\in\Delta_{C_4}$. 
Suppose $w\in\W_{\A, B}$ 
is compatible with $(\la,\pi)$ and such that $row(A(w))<C_5$. 
Then 
$$
\frac 1{C_6}\leq \frac{{\bf m}(C(w))}{\Prob(w|(\la,\pi))}\leq C_6
$$
\end{corollary}

\begin{corollary}
For any $C_7>0$, $C_8>0$ $C_9>0$, there exists $C_{10}>0$ such  that the following is true.

Suppose $(\la,\pi)\in\Delta_{C_7}$. 

Suppose $w\in\W_{\A, B}$ 
is compatible with $(\la,\pi)$ and furthermore satisfies 

$$
row(A(w))<C_8,\ \Delta(w)\subset \Delta_{C_9} 
$$

Then 
$$
\frac 1{C_{10}}\leq \frac{\Prob(C(w))}{\Prob(w|(\la,\pi))}\leq C_{10}
$$
\end{corollary}

\begin{corollary}
\label{condcomp}
Let $M$ be such that for any $n>M$ any two vertices in the Rauzy graph can be joined in 
$n$ steps.

Then for any $C_{17}>0$, $C_{18}>0$ $C_{19}>0$, there exists $C_{20}>0$ such  that the following is true.

Suppose $(\la,\pi)\in\Delta^+\cap \Delta_{C_{17}}$. 

Suppose $w\in W_{\A, B}^+$ satisfies

$$
row(A(w))<C_{18},\ \Delta(w)\subset \Delta^+\cap \Delta_{C_{19}} 
$$

Then for any $n\geq M$, we have
$$
\frac 1{C_{20}}\leq \frac{\Prob(C(w))}{\Prob^{(2n)}(w|(\la,\pi))}\leq C_{20}
$$
\end{corollary}

From the definition  (\ref{birkhoffmetric}) of the Hilbert metric it easily follows that for any 
$\la,\lap\in\Delta_{m-1}$ we have

\begin{equation}
\label{dist}
e^{-d(\la,\lap)}\lap_i\leq \la_i\leq e^{d(\la,\lap)}\lap_i.
\end{equation} 

\begin{proposition}
\label{lip1} 
Assume $\la,\lap\in\Delta_{\pi}^+$. 
Then
$$
\exp(-md(\la,\lap))\leq \frac{\rho(\la, \pi)}{\rho(\lap, \pi)}\leq \exp(md(\la,\lap))
$$
\end{proposition}

Proof. Indeed, there exist linear forms 
$$
l_i^{(j)}(\la)=\sum_{k=1}^m a^{(j)}_{ik}\la_k,
$$ 

where $a^{(j)}_{ik}$ are nonnegative integers (in fact, 
either $0$ or $1$, but we do not need this here),

such that 

$$
\rho(\la,\pi)=\sum_{j=1}^s \frac 
1{l_1^{(j)}(\la)l_2^{(j)}(\la)\dots l_m^{(j)}(\la)}.
$$

Clearly, if for all $i=1, \dots,m$ and some $\alpha>0$, we have
$\alpha^{-1}\la_i\leq\lap_i\leq\alpha\la_i$, then 

$$
\alpha^{-m}\leq \frac{\rho(\la,\pi)}{\rho(\lap, \pi)}\leq \alpha^m,
$$
and the Proposition is proved.

For similar reasons we have 

\begin{proposition}
\label{lip2} 
Assume $\la,\lap\in\Delta_{\pi}^+$ and let $A$ be an arbitrary matrix with nonnegative integer entries. 
Then
$$
\exp(-md(\la,\lap))\leq \frac{\rho(A\la, \pi)}{\rho(A\lap, \pi)}\leq \exp(md(\la,\lap))
$$
\end{proposition}

From these propositions and the formula \ref{trpr} we obtain

\begin{corollary}
\label{lip3}
Let $c\in\A$ be compatible with $\pi$. Then for any $\la,\lap\in\Delta^+_{\pi}$ we have

$$
\exp(-2md(\la,\lap))\leq \frac{\Prob(c|(\la,\pi))}{\Prob(c|(\lap,\pi))}\leq \exp(2md(\la,\lap))
$$

\end{corollary}

This Corollary implies the following

\begin{lemma}
\label{diam}
Let $w\in W_{\A,B}^+$ be such that the cylinder $C(w)$ has finite Hilbert diameter. 

Then for any $c$ compatible with $w$ and any $(\la_0,\pi)\in C(w)$ we have

$$
\exp (-2m\ diam C(w)) \leq \frac{\Prob(c|(\la_0,\pi))}{\Prob(\omega_0=c|\omega|_{[1, |w|]}=w)}
\leq \exp(2m\ diam C(w))
$$

\end{lemma}

Proof: We  have

$$
\nu(C(cw))=\int_{C(w)} \Prob(c|(\la,\pi))d\nu(\la,\pi)
$$

Let $d=diam C(w)$.
For any $(\la_,\pi), (\lap, \pi)\in C(w)$, we have, by Corollary \ref{lip3},  

$$
\exp(-2md)\leq \frac{\Prob(c|(\la,\pi))}{\Prob(c|(\lap,\pi))}\leq \exp(2md).
$$
 
Fix an arbitrary $(\la_0,\pi)\in \Delta_w$.

Then, from the above, 

$$
\nu(C(w))P(c|(\la_0,\pi))\exp (-2md)\leq 
\int_{C(w)} P(c|(\la,\pi))d\nu(\la,\pi)\leq 
$$
$$
\leq\nu(C(w))P(c|(\la_0,\pi))\exp (2md),
$$

and, since, by definition, we have  
$$
\Prob(\omega_0=c|\omega|_{[1, |w|]}=w)=\frac{\Prob(cw)}{\Prob(w)},
$$
the Lemma is proved.

For $N\in {\mathbb N}$ and $A\subset \Delta(\R)$, we denote
$\Prob^{(N)}(A|(\la,\pi))=\Prob((\la(-N), \pi(-N))\in A|(\la(0), \pi(0))=(\la,\pi))$; 
for $w\in\wab$, we write $\Prob^{(N)}(w|(\la,\pi))=\Prob^{(N)}(\Delta(w)|(\la,\pi))$.

\begin{lemma}
\label{compact} 
Let $M$ be a number such that for any 
$N\geq M$ any two vertices of the Rauzy graph can be connected in $N$ steps.
For any $\gamma>0$,  $N\geq M$  there exists a constant $C_0$ depending only on $\gamma$ and $N$
such that for any word $w\in\W_{\A,B}^+$ and any $(\la,\pi)\in\Delta_{\gamma}$
$$
\Prob^{(2N)}(w|(\la,\pi))\geq \frac{C_0}{|A(w)\la|^m}
$$
\end{lemma}

Proof: 

Let $w=w_1\dots w_{2n}$, and let 
$w_{2n}=(a, m_1, \pi_1)$.

Let $\pip_1\pip_2\dots \pip_{2N}$ 
a path of length $2N$ between $\pi$ and $\pi_1$ 
(here $\pip_1=\pi$, $\pip_{2n}=\pi_1$, $\pi_{2k+1}=a\pi_{2k}$, $\pi_{2k+2}=b\pi_{2k+1}$.

Denote $w_{n+2i+1}=(a, 1,\pi_{2i+1})$, $w_{n+2i}=(b,1,\pi_{2i})$. 
In other words, the word $=w_{2n+1}\dots w_{2n+2N}\in\W_{\A,B}$ is the word correspoding to the 
path $\pip_1\pip_2\dots \pip_{2N}$ in the Rauzy graph.
Then $w^{\prime}=w_1\dots w_{2n+2N}$ is a word compatible with $(\la,\pi)$.
Besides, 
$$
|A(w_{2n+1}c_{n+2}\dots w_{2n+2N})|<(2N)^{(2N)}.
$$

We have

$$
P^{(2n)}(w|(\la,\pi))\geq 
P(w^{\prime}|(\la,\pi))=
\frac{\rho(T_{w^{\prime}}(\la), w^{\prime}\pi)}
{|A(w^{\prime})\la|^m\rho(\la,\pi)},
$$

There exists a universal constant $C_1$ such that 
$\rho(\lap,\pip)>C_1$ for any $(\lap,\pip)\in\Delta^+$ (the density of the invariant measure is bounded from 
below).

Then, $|A(w^{\prime})\la|^m\leq |A(w^{\prime})|^m\leq (2N)^{2mN}|A(w)|^m.$

Finally, there exists a $C_2$ depending on $c$ only such that
if $\la_i>c$ for all $i$ then $\rho(\la,\pi)>C_2$.

Combining all of the above, we obtain the result of the Lemma.

\section{Markov approximation and the Doeblin condition}

\subsection{Good cylinders}

Let ${\bf q}=q_1\dots q_l$ be a word such that all entries of the matrix 
$A({\bf q})$ are positive. Fix $\epsilon>0$ and let et $k_0$ be such that 

\begin{equation}
\label{mixing}
\Prob(\Delta({\bf q})\cap \gz^{-2n}\Delta({\bf q}))\geq \epsilon \ {\rm for}\ n>k_0.
\end{equation}

Note that, due to mixing, Corollary \ref{qlogsteps} implies the following

\begin{proposition}
\label{qmixing}
Let $\q\in W_{\A,B}$, $\q=q_1\dots q_l$ be such that 
all entries of the matrix $A(\q)$ are positive and that $\Delta(\q)\subset \Delta^+$.
Then there exist positive constants $K(\q),p(\q)$ 
such that the following is true for any $\epsilon>0$. 
Suppose $(\la, \pi)\in\Delta_{\epsilon}\cap \Delta^+$ and set $n$ to be the integer 
part of $K(\q)|\log \epsilon|$. Then  

$$
\Prob \{ (\la(-2n),\pi(-2n))\in\Delta(\q)|
(\la(0), \pi(0))=(\la,\pi))\}\geq p(\q).
$$
\end{proposition}

Take $k\geq k_0$. Let $r=2(K+1)k+2M$, where $K$ is the constant from the Lemma \ref{logsteps}
and $M$ is the connecting constant of the Rauzy graph from Lemma \ref{compact}. 

Let $\theta$, $0<\theta<1$ be arbitrary.
A word $w=w_1\dots w_k$ is called {\it good} if 
\begin{enumerate}
\item $\Delta(w)\subset \Delta_{\exp(-k)}$.
\item the word ${\bf q}$ appears at least $\frac{k^{\theta}}{l}$
times in $w$ (we only count disjoint appearances).
\end{enumerate} 
A word $w_1\dots w_r$ is called good if $w_1\dots w_k$ is good, 
a word $w_1\dots w_{Nr}$ is called good if all words 
$w_1\dots w_r$, $w_{r+1}\dots w_{2r}$, \dots $w_{(N-1)r+1}\dots w_{Nr}$
are good, and a word $w_1\dots w_{Nr+L}$, $L<r$, is good if 
$w_1\dots w_{Nr}$ is good and either $L<k$ or $w_{Nr+1}\dots w_{Nr+k}$ is good.   

We denote by $\G(N)$ the set of all good words of length $N$.

Let $$
\Delta(\G(N))=\cup_{w\in\G(N)}\Delta(w), 
$$
and 
$$
\Delta(B(N))=\Delta^+\setminus \Delta(\G(N))
$$

By Corollary \ref{sqroot}, there exist constants $C_{31}, C_{32}$ such that for all $r$ we have 

\begin{equation}
\label{measbad}
\Prob(\Delta(B(N))\leq C_{31}N\exp(-C_{32}r^{(1-\theta)/2}).
\end{equation}

and, for any $(\la,\pi)\in \Delta(\q)$, also 

\begin{equation}
\label{condmeasbad}
\Prob((\la(-1), \pi(-1))\in\Delta(B(N))|(\la(0), \pi(0)=(\la,\pi))\leq C_{31}N\exp(-C_{32}r^{(1-\theta)/2}).
\end{equation}

\subsection{Preliminary estimates for the Doeblin condition.}

From Corollary \ref{lip3} we deduce that there exists a  constant $C_{33}$ such that for any 
$(\la,\pi), (\la^{\prime}, \pi)\in\Delta({\bf q})$, and any word $w$ 
compatible with ${\bf q}$, we have 
$$
\frac1C_{33}\leq \frac{\Prob(w|(\la,\pi))}{\Prob(w|(\la^{\prime},\pi))}\leq C_{33}.
$$

Finally, by Lemma \ref{compact}, there exists a constant $C_{34}$ such that
for any $w\in\W_{\A, B}$ and for any $N>M$ we have 
$$
\frac1C_{34}\leq \frac{\Prob^{(2N)}(w|(\la,\pi))}
{\Prob^{(2N)}(w|(\la^{\prime},\pi))}\leq C_{34}.
$$

Take an arbitrary point $(\la,\pi)\in\Delta_{\q}$.
Define a new measure $\varphi$ on $\Delta^+$.
Namely, for a set $A\subset\Delta^+$ put  
\begin{equation}
\label{doebmeas}
\varphi(A)=\Prob(\la(-2M), \pi(-2M))\in A|\la(0), \pi(0)=(\la,\pi))
\end{equation}

\begin{lemma}
\label{predoeblin}
There exists a constant $\alpha>0$ such that the following is true for any $r$.
Let ${\cal C}_1, {\cal C}_2\in\G(r)$.

Then 
$$
\Prob({\bf \omega}|_{[1,r]}={\cal C}_1, \omega|_{[r+1, 2r]}\in\G(r)\ |{\bf \omega}|_{[2r+1,3r]}={\cal C}_2)\geq 
\alpha\varphi({\cal C}_1)
$$
\end{lemma}

Indeed, we have the following propositions:

\begin{proposition}
\label{gettoq}
There exist a constant $p_1$ such that the following is true for all $r$ and all $n\geq r$.

Let $C_2\in\G(r)$, $(\la,\pi)\in {\cal C}_2$.
Then 
$$
\Prob((\la(-2n), \pi(-2n))\in\Delta(\q)|(\la(0), \pi(0))=(\la,\pi))\geq p_1.
$$
\end{proposition}

This follows from the definition of a good cylinder and Corollary \ref{qlogsteps}.

\begin{proposition}

There exists a constant $p_2$ such that the following is true for all $k$. 
$$
\Prob(\omega|_{[1,r]}\in\G(r), \omega|_{[2M+1, l+2M+1]}=\q\ | \omega|_{[r+1, r+l+1]}=\q)\geq p_2 
$$
\end{proposition}

This follows from the estimates (\ref{measbad}),(\ref{condmeasbad}) on the measure of bad cylinbders and from Proposition 
\ref{qmixing}. 

\begin{proposition}
There exists a constant $p_3$ such that the following is true for all $r$.
Let $c_1\dots c_n\dots \in \Delta(\q)$.

$$
\Prob(\omega|_{[1,r]}=\C_1|\omega_{r+2M+1}=c_{1}, \omega_{r+2M+2}=c_{2}, \dots)\geq p_3\varphi(C_1)
$$
\end{proposition}

This follows directly from Lemma \ref{compact}.

The three Propositions imply  Lemma \ref{predoeblin}.

\subsection{Approximation by a Markov measure}

We define a new measure $\prt$ on the set $\G(r^2)$ of good cylinders of length $r^2$. 

Let $\C=c_1\dots c_{r^2}$ be a $(r, \theta)$-good cylinder. 
Set $\C_i=c_{ir+1}\dots c_{(i+1)r}$.

Define 
$$
\prt(\cc)=\Prob(\omega|_{[1,r]}=\C_1|\omega|_{[r+1, 2r]}=\C_2) 
\Prob(\omega|_{[r+1, 2r]}=\C_2|\omega|_{[2r+1, 3r]}=\C_3)\dots 
\Prob(\omega|_{[r^2-r+1, r^2]}=\C_r).
$$

If $D$ is not a good cylinder, then $\prt(D)=0$. 

Normalize to get a probability measure:

$$
\Prt(\C)=\frac{\prt(\C)}{\sum_{\D\in\G(r^2)}\prt(\D)}.
$$

$\Prt$ is a Markov measure of memory $r$ (in general, 
non-homogeneous), as is shown by 
the following well-known Lemma \cite{bunimsinai}. 
\begin{lemma}
For any $k$, $0<k<r$, we have 
$$
\Prt(\omega|_{[kr+1, (k+1)r]}=\C_k|\omega|_{[(k+1)r+1, r^2]})=\C_{k+1}\dots \C_r)=
$$
$$
\Prt(\omega|_{[kr+1, (k+1)r]}=\C_k|\omega|_{[(k+1)r+1, (k+2)r]})=\C_{k+1}).
$$
\end{lemma}

From the H{\"o}lder property for the transition probability, we have

\begin{proposition}

There exist constants $C_{41}, C_{42}$ such that the following is true for any $r$. 

Let $c_1\dots c_n\dots \in\Omega_{\A, B}$ and assume 
$c_{n+1}\dots c_{n+r}\in\G(r)$. Then 

$$
\exp(-C_{41}\exp(-C_{42}k^{\theta}))\leq 
$$
$$\leq\frac{P(\omega_1=c_1, \dots, \omega_n=c_n|\omega_{n+1}=c_{n+1}, \dots, \omega_{n+r}=c_{n+r})}
{\Prob(\omega_1=c_1, \dots, \omega_n=c_n|\omega_{n+1}=c_{n+1}, \dots, \omega_{n+i}=c_{n+i}, \dots)}  
\leq 
$$
$$
\leq \exp(C_{41}\exp(-C_{42}k^{\theta}))
$$
\end{proposition}

\begin{corollary}
There exist constants $C_{43}, C_{44}$ such that the following is true for any $r$.
Let $A\in\F_n$, let $c_{n+1}\dots c_{n+i}\dots\in\Omega_{\A}$,  
and assume $c_{n+1}\dots
c_{n+r}\in\G(r)$. Then

$$
\exp(-C_{43}\exp(-C_{44}k^{\theta}))\leq \frac
{\Prob(A|\omega_{n+1}=c_{n+1}, \dots,
\omega_{n+r}=c_{n+r})}
{\Prob(A|\omega_{n+1}=c_{n+1}, \dots,
\omega_{n+i}=c_{n+i}, \dots)}
\leq \exp(C_{43}\exp(-C_{44}k^{\theta}))
$$
\end{corollary}

Applying $l$ times, we obtain

\begin{lemma}
There exist constants $C_{45}, C_{46}, C_{47}, C_{48}$ such that the following is true for any $r$.
Let $c_1\dots c_{r^2} \in\G(r^2)$. Then for any $l$, $1\leq l\leq r$, 
we have
$$
\exp(-C_{45}l\exp(-C_{46}k^{\theta}))\leq 
$$
$$
\leq\frac
{\Prob(\omega_1=c_1, \dots, \omega_{lr}=c_{lr}|\omega_{lr+1}=c_{lr+1}, \dots,
\omega_{r^2}=c_{r^2})}
{\prt (\omega_1=c_1, \dots, \omega_{lr}=c_{lr}|\omega_{lr+1}=c_{lr+1}, 
\dots,\omega_{r^2}=c_{r^2})}
\leq 
$$
$$
\leq \exp(C_{45}l\exp(-C_{46}k^{\theta}))
$$

and 

$$
\exp(-C_{47}l\exp(-C_{48}k^{\theta}))\leq \frac{\Prob(\omega_1=c_1, \dots, \omega_{lr}=c_{lr})}
{\prt (\omega_1=c_1, \dots, \omega_{lr}=c_{lr})}
\leq \exp(C_{47}l\exp(-C_{48}k^{\theta}))
$$
\end{lemma}

Summing over cylinders of length $lr$, we obtain

\begin{corollary}
There exist constants $C_{49}, C_{50}$ such that the following is true for any $r$.
Let $c_1\dots c_{r^2} \in\G(r^2)$. Then for any $l$, $1\leq l\leq r$,
and any $A\in\F_{lr}$, we have 
$$
\exp(-C_{49}l\exp(-C_{50}k^{\theta}))\leq \frac
{\Prob(A\cap \G(lr)|\omega_{lr+1}=c_{lr+1}, \dots,
\omega_{r^2}=c_{r^2})}{\prt (A|\omega_{lr+1}=c_{lr+1},
\dots,\omega_{r^2}=c_{r^2})}
\leq \exp(C_{49}l\exp(-C_{50}k^{\theta}))
$$

and

$$
\exp(-C_{49}l\exp(-C_{50}k^{\theta}))\leq \frac
{\Prob(A\cap \G(lr))}
{\prt (A)}
\leq \exp(C_{49}l\exp(-C_{50}k^{\theta}))
$$

\end{corollary}

Using (\ref{measbad}), we can estimate the total mass of the measure $\prt$.

\begin{corollary}
There exist constants $C_{51}, C_{52}$ such that for any $r$ we have
$$
\prt(\G(r^2))\geq \exp(-C_{51}r\exp(-C_{52}k^{(1-\theta)/2}))
$$
\end{corollary}

We now have normalized versions of previous statements.

\begin{corollary}
There exist constants $C_{53}, C_{54}, C_{55}, C_{56}$ such that the following is true for any $r$.
Let $c_1\dots c_{r^2} \in\G(r^2)$. Then for any $l$, $1\leq l\leq r$,
and any $A\in\F_{lr}$, we have
$$
\exp(-C_{53}l\exp(-C_{54}k^{\theta})-C_{55}r\exp(-C_{56}k^{(1-\theta)/2}))\leq 
$$
$$\leq \frac
{\Prob(A\cap \G(lr)|\omega_{lr+1}=c_{lr+1}, \dots,  
\omega_{r^2}=c_{r^2})}{\Prt (A|\omega_{lr+1}=c_{lr+1},
\dots,\omega_{r^2}=c_{r^2})}
\leq 
$$
$$\leq \exp(C_{53}l\exp(-C_{54}k^{\theta})+C_{55}r\exp(-C_{56}k^{(1-\theta)/2})
$$

and

$$
\exp(-C_{53}l\exp(-C_{54}k^{\theta})-C_{55}r\exp(-C_{56}k^{(1-\theta)/2}))\leq 
$$
$$\leq\frac
{\Prob(A\cap \G(lr))}{\Prt (A)}
\leq 
$$
$$
\leq\exp(C_{53}l\exp(-C_{54}k^{\theta})+C_{55}r\exp(-C_{56}k^{(1-\theta)/2}).
$$
\end{corollary}

Using the Markov approximation, we can  estimate conditional measure of good cylinders for the measure ${\mathbb P}$:

\begin{corollary}
There exist constants $C_{57}, C_{58}, C_{59}, C_{60}$ such that the following is true for any $r$.
Let $c_1\dots c_{r^2} \in\G(r^2)$. Then for any $l$, $1\leq l\leq r$,
we have
$$
\Prob((\omega_1\dots \omega_{lr})\in\G(lr)|\omega_{lr+1}=c_{lr+1}, 
\dots,\omega_{r^2}=c_{r^2})\geq 
\exp(-C_{57}l\exp(-C_{58}k^{\theta})-C_{59}r\exp(-C_{60}k^{(1-\theta)/2}))
$$
\end{corollary}

Proof: Indeed, 
$$
\Prt((\omega_1\dots \omega_{lr})\in\G(lr)|\omega_{lr+1}=c_{lr+1},
\dots,\omega_{r^2}=c_{r^2})=1.
$$

\subsection{Doeblin Condition}

\begin{proposition}
There exists $C_{61}$ such that the following holds for any $r$.
For any  $\C_1\subset \Delta({\q})$, $C_2\subset \Delta_{{\q}}$, and any 
$\C_3\in \G(r)$, we have  either 

$$
\frac 1{C_{61}}\leq \frac{\prt(C_3|C_2)}{\prt(C_3|C_1)}\leq C_{61},
$$

or $\prt(C_3|C_2)=\prt(C_3|C_1)=0$. 
\end{proposition}

Considering $n$-step transition probabilities, we obtain

\begin{proposition}
\label{boundratio}
There exists a constant $C_{62}$ such that the following holds for any $r$.
For any  $\C_1\subset \Delta({\q})$, $C_2\subset \Delta_{{\q}}$ any 
$\C_3\in \G(r)$, and any $n\geq M$, we have

$$
\frac 1{C_{62}}\leq \frac
{\prt(\omega|_{[1,r]}=C_1|\omega_{[2n+r,2n+2r]}=C_2)}
{\prt(\omega|_{[1,r]}=C_1|\omega_{[2n+r,2n+2r]}=C_3)}\leq C_{62}.
$$

\end{proposition}

Now, mixing, Proposition \ref{qmixing} and Proposition \ref{gettoq},
and the definition oif a good cylinder 
imply that 

\begin{proposition}
There exists a constant $C_{63}$ such that the following holds for any $r$.
For any $\C_1, \C_2, \C_3\in \G(r)$ we have
$$
\frac 1{C_{63}}\leq \frac
{\prt(\omega|_{[1,r]}=C_1|\omega_{[2r,3r]}=C_2)}
{\prt(\omega|_{[1,r]}=C_1|\omega_{[2r,3r]}=C_3)}\leq C_{63}.
$$
\end{proposition}

Now let $c_1\dots c_{r^2} \in\G(r^2)$. Denote $\C_i=c_{ir+1}\dots c_{(i+1)r}$.
Lemma \ref{predoeblin}, together with the above estimates, implies the following

\begin{corollary}
\label{doeblin}
There exist constants $C_{71}, C_{72}$ such that the following is true.
For any $l$, $1\leq l\leq r$,
we have
$$
\Prob(\omega|_{[1,lr]}\in\G(lr), 
\omega|_{[lr+1, (l+1)r]}=\C_l, 
\omega|_{[(l+1)r+1,(l+2)r]}\in\G(r)|
\omega_{(l+2)r+1, (l+3)r])}=\C_3)\geq 
C_{71}\times \varphi(\C_l) 
$$

and 

$$
\Prt(\omega|_{[1,lr]}\in\G(lr), 
\omega|_{[lr+1, (l+1)r]}=\C_l, 
\omega|_{[(l+1)r+1,(l+2)r]}\in\G(r)|
\omega_{(l+2)r+1, (l+3)r])}=\C_3)\geq C_{72}\times \varphi(\C_l) 
$$
\end{corollary}

This is the Doeblin Condition for the measure $\Prt$ (see \cite{sinai}, \cite{bunimsinai}, \cite{doob}). 
The Doeblin Condition implies that there exist constants $C_{73}, C_{74}$ such that 
for any $\C_1, \C_2\in \G(r)$, we have

$$
\exp(-C_{73}\exp(-C_{74}r))\leq 
\frac{\Prt(\omega|_{[1,r]}=\C_1|\omega|_{[r^2, r^2+r]}=\C_2)}{\Prt(\C_1)}\leq \exp(C_{73}\exp(-C_{74}r)),
$$

whence we obtain

\begin{proposition}
\label{correl}
There exist constants $C_{81}, C_{82}, C_{83}, C_{84}$ such that the following is true for any $r$.
$$
\exp(-C_{81}(\exp(-C_{82}r)+\exp(-C_{83}r^{\theta})+\exp(-C_{84}r^{(1-\theta)/2})))\leq 
$$
$$
\leq\frac{\Prob(\omega|_{[1,r]}=\C_1|\omega|_{[r+1, r^2]}\in \G(r^2-r), \omega|_{[r^2, r^2+r]}=\C_2)}{\Prob(\C_1)}\leq
$$
$$
\leq 
\exp(C_{81}\exp(-C_{82}r)+\exp(-C_{83}r^{\theta})+\exp(-C_{84}r^{(1-\theta)/2)}))).
$$
\end{proposition}

Moreover, in view of mixing, Proposition \ref{qmixing}, 
and Proposition \ref{gettoq}, the same estimate, upto a constant, takes place for any $n\geq r^2$.

\begin{proposition}
\label{correln}
There exist constants $C_{85}, C_{86}, C_{87}, C_{88}$ such that the following is true for all $r$ and all $n\geq r^2$.
$$
\exp(-C_{85}(\exp(-C_{86}r)+\exp(-C_{87}r^{\theta})+\exp(-C_{88}r^{(1-\theta)/2}))))\leq 
$$
$$
\leq\frac{\Prob(\omega|_{[1,r]}=\C_1|\omega|_{[r+1, n]}\in \G(n-r), \omega|_{[n, n+r]}=\C_2)}{\Prob(\C_1)}\leq
$$
$$
\leq 
\exp(C_{85}(\exp(-C_{86}r)+\exp(-C_{87}r^{\theta})+\exp(-C_{88}r^{(1-\theta)/2)})))).
$$
\end{proposition}

\section{Approximation of H{\"o}lder Functions and Completion of the Proof of Theorems \ref{mainresult}, \ref{zipdecay}, \ref{cltpt}.}

We shall prove the decay of correlations for 
a slightly more general class of functions on $\Delta(\R)$ than H{\"o}lder functions. 
(we shall need this slightly more general class in the proof of the Central Limit Theorem).

Namely, we shall only require that a function be H{\"o}lder in restriction to 
cylinders of some given length and we shall also allow a moderate growth 
of the H{\"o}lder constant at infinity.

Formally, say that a function $\phi:\Delta(\R)\to{\mathbb R}$ is {\it weakly $l, \alpha$-H{\"o}lder} if 
the following holds. Let $k$ be a positive integer, and let  
$w\in\wab$, $|w|\leq l$ be such that $\Delta(w)\subset \Delta_{\exp(-k)}$. 
Then there exists a constant $C(\phi)$ such that for any $(\la,\pi), (\lap, \pi)\in\Delta(w)$,
we have 
$$
|\phi(\la,\pi)-\phi(\lap,\pi)|\leq C k d(\la,\lap)^{\alpha}.
$$
 
The smallest such $C$ for a given $\phi$ will be denoted $C_{l,\alpha}^{weak}(\phi)$. 
Clearly, if $\phi$ is H{\"o}lder with exponent $\alpha$, 
then it is also  weakly $l, \alpha$-H{\"o}lder for any $l$ and 
 $C_{l,\alpha}^{weak}(\phi)\leq C_{\alpha}(\phi)$. 

Recall that  ${\cal B}_n$ is the $\sigma$-algebra of sets of the form $\gz^{-n}(A)$, $A\subset \Delta(\R)$.

To prove the decay of correlations, it suffices to estimate the $L_2$-norm of 
$E(\phi|{\cal B}_{2n})$ for a given {\it weakly $l$-$\alpha$-H{\"o}lder} $\phi$.  

It will be convenient to assume that $\phi\geq 1$ (by linearity, it suffices to 
consider that case). 

\begin{proposition}
\label{condexpest}

Let $\theta\in{\mathbb R}$, $0<\theta<1$. Let $p>2$ and $\alpha>0$. 
There exist constants $C_{91}, C_{92}$, $C_{93}$ 
such that the following is true for any $r$ and any $n\geq r^2$.

Let $l\leq r$. Let $\phi\in L_p(\Delta(\R)^+, \nu)$ be weakly $l, \alpha$-H{\"o}lder and satisfy $\phi\geq 1$. 

Then $\phi=\phi_1+\phi_2+\phi_3$ where

\begin{enumerate} 

\item $\phi_1\geq 1$ on $\G(n)$ and $\phi_1=\phi_2=0$ on $\Delta(B(n))$.

\item for any $(\la,\pi)\in \G(n)$, we have 

$$|\frac{E(\phi_1|\F_n)(\la,\pi)}{E(\phi_1)}-1|\leq\exp(-C_{91}(r^{(1-\theta)/2}+r^{\theta}).$$ 

\item for $(\la,\pi)\in G(n)$, we have $|\phi_2|\leq C_{l,\alpha}^{weak}(\phi)\exp(-C_{92}r^{\theta}).$

\item $||\phi_3||_{L_2}\leq \exp(-C_{93}r^{(1-\theta)/2})||\phi||_{L_p}.$
\end{enumerate}

\end{proposition}

Proof: For any good word $w=w_1\dots w_{n+r}$, consider its beginning $w_1\dots w_r$ and choose 
a point $x_{w_1\dots w_r}\in \Delta(w_1\dots w_r)$. 

Denote by $\chi_{\Delta(w)}$ the characteristic function of $\Delta(w)$ and set

$$
{\phi_1}=\sum_{w\in \G(n+r)}\phi(x_{w_1\dots w_r})\chi_{\Delta(w)}. 
$$

Proposition \ref{correln}  yields the required properties of $\phi_1$ 
(note that we sum over all good words of length $n+r$ in order to be able to apply the Proposition).

We set $\phi_2=(\phi-\phi_1)\chi_{G(n+r)}$ and 
$\phi_3=\phi\chi_{\Delta(B(n+r))}$.
The estimate for $\phi_2$ is satisfied by the definition of a H{\"o}lder function. 

Finally, we have   
$$
||\phi_3||_{L_2}^2=
E(|\phi\chi_{\Delta(B(n))}|^2),
$$

whence, by H{\"o}lder's inequality, using the estimate (\ref{measbad}), 
we obtain the desired estimate for $\phi_3$, and  the Proposition is proved completely.

Proposition \ref{condexpest} with $\theta=1/3$ yields Theorem \ref{mainresult}.

We now complete the proof of Theorem \ref{zipdecay}.

For a word $w\in W_{\A,B}$, $|w|=2n+1$, $w=w_1\dots w_{2n+1}$, denote
$C^{[-n,n]}(w)=\{\omega\in\Omega_{\A,B}^{\mathbb Z}:\{\omega_{-n}=w_{1}, \dots,  
\omega_n=w_{2n+1}\}$ and set ${\overline \Delta}(w)={\overline \Phi}^{-1}C^{[-n,n]}(w)$. 
Denote by ${\cal B}_{[-n,n]}$ the sigma-algebra 
generated by ${\overline \Delta}(w)$ for all $w\in W_{\A,B}$.

Also, for $\epsilon>0$, denote 
$$
{\overline \Delta}_{\epsilon}=\{(\la,h,a,\pi)\in {\overline \Delta}(\R): \la\in\Delta_{\epsilon}.
$$

Again, we shall prove the Theorem for a slightly larger class of functions. 

Say that a function $\phi:{\overline \Delta}(\R)\to{\mathbb R}$ is {\it weakly $l, \alpha$-H{\"o}lder} if 
the following holds. Let $k$ be a positive integer, and let  
$w\in\wab$, $|w|\leq 2l+1$ be such that ${\overline \Delta}(w)\subset \Delta_{\exp(-k)}$. 
Then there exists a constant $C(\phi)$ such that for any $(\la,h,a,\pi), (\lap,h^{\prime}, a^{\prime}, 
\pi)\in{\overline \Delta}(w)$,
we have 
$$
|\phi(\la,h,a, \pi)-\phi(\lap, ,h^{\prime}, a^{\prime}, \pi)|\leq C k 
d((\la,h,a, \pi),(\lap, ,h^{\prime}, a^{\prime}, \pi))^{\alpha}.
$$
 
The smallest such $C$ for a given $\phi$ will be denoted $C_{l,\alpha}^{weak}(\phi)$. 
Clearly, if $\phi$ is H{\"o}lder with exponent $\alpha$, 
then it is also  {\it weakly $l, \alpha$-H{\"o}lder} for any $l$ and 
 $C_{l,\alpha}^{weak}(\phi)\leq C_{\alpha}(\phi)$.

Denote  by ${\overline G}(2n+1)$ the union of all 
${\overline \Delta}(w)$ for  good $w$, by ${\overline B}(2n+1)$
the  complement of  ${\overline G}(2n+1)$. 

\begin{proposition}
\label{approxcyl}

Let $\theta\in{\mathbb R}$, $0<\theta<1$. Let $p>2$ and $\alpha>0$. 
There exist constants $C_{101}, C_{102}$, 
such that the following is true for any $r$ and any $n\geq r^2$.

Let $l\leq r$. Let $\phi\in L_p({\overline \Delta}(\R)^, {\overline \nu})$ be weakly 
$l,\alpha$-H{\"o}lder and satisfy $\phi\geq 1$. 
Then there exist functions $\phi_1$, $\phi_2$, $\phi_3$ such that

\begin{enumerate}
\item $\phi=\phi_1+\phi_2+\phi_3$.
\item $\phi_1$ is ${\cal B}_{[-n,n]}$-measurable and supported on ${\overline G}(2n+1)$. 
\item $|\phi_2|\leq C_{101}C_{\alpha}(\phi)\exp(-r^{(1-\theta)/2}+r^{\theta})$.
\item $|\phi_3|_{L_2}\leq C_{102}\exp(-r^{(1-\theta)/2)}||\phi||_{L_p}$.
\end{enumerate}
\end{proposition}

For any good $w$, $|w|=2n+1$, take an arbitrary point $x_w$ in 
${\overline \Delta}(w)$.
Set
$$
\phi_1=\sum_{w\in G(2n+1)}\phi(x_{w})\chi_{{\overline \Delta}(w)}, 
$$
$$
\phi_2=(\phi-\phi_1)\cdot \chi_{G(2n+1)},
$$
$$
\phi_3=\phi\cdot \chi_{B(2n+1)},
$$
and the Proposition is proved.

Proposition \ref{approxcyl} with $\theta=1/3$ yields Theorem \ref{zipdecay}.

It remains to establish the Central Limit Theorem for the flow $P^t$.
Consider the special function ${\tilde \tau}$ of the flow $P^t$ over the transformation 
$\F$. Note that ${\tilde \tau}(\la,h,a,\pi)$ only depends on $(\la,\pi)$.
Consider the restriction of ${\tilde \tau}$ on a cylinder of the form $\Delta(w_1)$,  
$w_1\in\A$. Then there exist distinct $j(1), \dots, j(l)\in\{1, \dots, m\}$ such that

$$
{\tilde \tau}(\la, \pi)=\log(\la_{j(1)}+\la_{(j(2)}+\dots+\la_{j(l)}),
$$

which shows that the function ${\tilde \tau}$, restricted to an arbitrary $\Delta(w_1)$
is Lipshitz with respect to the Hilbert metric on $\Delta(\R)$.

Now for a H{\"o}lder $\phi$ consider the function

$$
{\tilde \phi}(x)=\int_0^{{\tilde \tau(x)}}\phi(P^tx).
$$

For any $k>1$, if $(\la,\pi)\in \Delta_{\exp(-k)}$, then, by definition, ${\tilde \tau}(\la, \pi)\leq k$. 
Therefore, if $\phi$ is H{\"o}lder of exponent $\alpha$, then ${\tilde \phi}$ is weakly $1, \alpha$-H{\"o}lder.

It is easy to see that ${\tilde \tau}(\la,\pi)\in L_r(\Delta(\R),\nu)$ for any $r>1$, whence,
if $\phi\in L_p(\Omega_0(\R), \mu_{\R})$ for some $p>2$, then there exists $p^{\prime}>2$
such that the function
$$
{\tilde \phi}(x)=\int_0^{{\tilde \tau(x)}}\phi(P^tx)
$$

satisfies ${\tilde \phi}\in L_{p^{\prime}}(\Y^{\pm}, {\overline \nu})$.

Therefore, the Theorem of Melbourne and T{\"o}r{\"o}k \cite{torok} implies Theorem \ref{cltpt}, the Central Limit 
Theorem for the flow $P^t$.

{\bf Acknowledgements.}

I am deeply grateful to  Yakov G. Sinai, who 
introduced me to interval exchange transformations,
explained to me the method of Markov approximations, and 
encouraged me in every possible way
as the work progressed (more importantly, when it did not).

I am deeply grateful to Alexander Eskin, who suggested to me  
the problem of the decay of correlations for the induction map.

I am deeply grateful to Giovanni Forni who introduced me to Teichm{\"u}ller theory.

I am deeply grateful to Corinna Ulcigrai and Pavel Batchourine for their 
suggestions that have been of invaluable help to me.

I am deeply grateful to Jayadev Athreya, Valdo Durrleman, Charles L. Fefferman,
Boris M. Gurevich, Carlangelo Liverani, Michael Ludkovski, Ian Melbourne, Klaus Schmidt, 
Andrei T{\"o}r{\"o}k and Anton V. Zorich for useful discussions.

Part of this work was done at The Erwin Schr{\"o}dinger Institute in 
Vienna, at The Institute of Mathematics "Guido Castelnuovo" of the 
University of Rome "La Sapienza", and at the CIRM-IML in Marseille. 
I am deeply grateful to these institutions for their hospitality.

\end{document}